\documentclass[12pt]{amsart}
\usepackage[margin=1in]{geometry}                
\usepackage{graphicx,color}
\usepackage{amssymb,amsmath,amsthm}
 \usepackage{amsaddr}
 
\newtheorem{conjecture}{Conjecture}
\newtheorem{proposition}{Proposition}
\newtheorem*{proposition*}{Proposition}

\newcommand{\beq}{\begin{equation}}
\newcommand{\eeq}{\end{equation}}

\newcommand{\ds}{\displaystyle}

\renewcommand{\vec}[1]{\boldsymbol{#1}}

\newcommand{\EE}{\mathbb{E}}
\newcommand{\RR}{\mathbb{R}}

\begin{document}

\title[Metastability of the Nonlinear Wave Equation]
  {Metastability of the Nonlinear Wave Equation: \\
  Insights from Transition State Theory}

\author{Katherine A Newhall}
\address{Department of Mathematics, University of North Carolina at Chapel Hill}
\author{Eric Vanden-Eijnden}
\address{Courant Institute of Mathematical Sciences,
New York University}

  \date{\today} 

\begin{abstract}

  This paper is concerned with the long-time dynamics of the nonlinear wave equation in one-space dimension,
$$
u_{tt} - \delta^2 u_{xx} +V'(u) =0 \qquad x\in [0,1]
$$
where $\delta>0$ is a parameter and $V(u)$ is a potential bounded from below and growing at least like $u^2$ as $|u|\to\infty$. Infinite energy solutions of this equation preserve a natural Gibbsian invariant measure and when the potential is double-welled, for example when $V(u) = \tfrac14(1-u^2)^2$, there is a regime such that two small disjoint sets in the system's phase-space concentrate most of the mass of this measure. This suggests that the solutions to the nonlinear wave equation can be metastable over these sets, in the sense that they spend long periods of time in these sets and only rarely transition between them. Here we quantify this phenomenon by calculating exactly via Transition State Theory (TST) the mean frequency at which the solutions of the nonlinear wave equation with initial conditions drawn from its invariant measure cross a dividing surface lying in between the metastable sets. We also investigate numerically how the mean TST frequency compares to the rate at which a typical solution crosses this dividing surface. These numerical results suggest that the dynamics of the nonlinear wave equation is ergodic and rapidly mixing with respect to the Gibbs invariant measure when the parameter $\delta$ in small enough. In this case successive transitions between the two regions are roughly uncorrelated and their dynamics can be coarse-grained to jumps in a two-state Markov chain whose rate can be deduced from the mean TST frequency. This is a regime in which the dynamics of the nonlinear wave equation displays a metastable behavior that is not fundamentally different from that observed in its stochastic counterpart in which random noise and damping terms are added to the equation. For larger $\delta$, however, the dynamics either stops being ergodic, or its mixing time becomes larger than the inverse of the TST frequency, indicating that successive transitions between the metastable sets are correlated and the coarse-graining to a Markov chain fails.

\end{abstract}

\maketitle

\section{Introduction}
\label{sec:intro}

Metastability is a feature commonly displayed by dynamical systems~\cite{bovierbook}. One of its simplest manifestation is in the context of finite-dimensional systems whose evolution is governed by the so-called overdamped Langevin equation, that is, steepest descent over a given energy subject to thermal fluctuations modeled by an additive white-noise driving term. In such systems, if the energy landscape possesses several local minima separated by barriers much higher than the available thermal energy, with high probability the solution will spend a long time in a basin associated with the energy minimum before making a transition to another such basin. The time scale of these transitions follows the Arrhenius law, i.e. it is exponentially large in the energy barrier the system must surmount to escape the basin, measured in units of temperature. This prediction can be confirmed rigorously using Freidlin-Wentzell large deviation theory~\cite{FWbook}, which also predicts the mechanism of the transitions. More refined estimates for the transition time that include its prefactor can be obtained using e.g.~the potential theoretic approach to metastability introduced by Bovier and collaborators~\cite{bovierbook,Bovier:2001dk,Bovier1,Bovier2}, and extended by many others~\cite{Barret,Huisinga:2004vr,EVEMaria2008}.

These results can be extended to infinite-dimensional dynamical systems, e.g. the stochastic Allen-Cahn equation~\cite{Berglund:2013fz,Faris:1982gs,EVEMaria2008}. They can also be extended to other types of dynamics, e.g.~if inertial effects are included. In this case, the system's evolution can be modeled by the (inertial) Langevin equation, that is, Hamilton's equations in which a linear damping term and a white-noise forcing term are added in the equation for the momenta. More complicated is the case of Hamiltonian systems in which there is no damping nor thermal noise, since the dynamics of such systems is conservative and deterministic. Still, metastability can arise in this context as well. This is most easily pictured in finite dimensional systems, assuming that the constant energy surface the system is evolving upon consists of two large regions connected by a narrow bottleneck: this arises e.g.~if the total energy of the system is only slightly higher than the potential energy barrier between two regions in position space.  If the dynamics is ergodic, the system can and will eventually cross this bottleneck, but in order to do so it must give up most of its kinetic energy to overcome the potential energy, which is highly unlikely to happen. This picture is intuitively appealing, though it may be hard to justify rigorously since ergodicity is typically very hard to prove. On top of this, it also becomes harder to use and justify as the dimensionality of the system under consideration increases. Indeed, energy typically is an extensive property that grows with the system's size, meaning that large systems always have much more energy that what is needed to cross a given potential energy barrier.  In particular, the picture above seems to become moot for infinite-dimensional Hamiltonian systems, since their total energy can be infinite.  The purpose of the present paper is to investigate metastability in the context of one such system.

Specifically, we will study the nonlinear one-dimensional wave equation for $u(x,t)$
\begin{equation}
  \label{full_pde00}
  u_{tt} - \delta^2 u_{xx} + V'(u) = 0
\end{equation}
for $x\in[0,1]$, with either Dirichlet or Neumann boundary conditions at $x=0$ and $x=1$, and $V(u)>Cu^2$ for some $C>0$ as $|u|\to\infty$. In \eqref{full_pde00}, $\delta$ measures the typical spatial scale over which $u(x,t)$ varies: this parameter can be absorbed in the domain size after rescaling via $\xi = x/\delta \in [0,1/\delta]$.  Defining $p(x,t) = u_t(x,t)$, we can write \eqref{full_pde00} as
\begin{equation}
  \label{full_pde0}
    u_t = p,\qquad
    p_t = \delta^2 u_{xx} - V'(u)
\end{equation}
which are Hamilton's equations associated with the infinite-dimensional Hamiltonian
\begin{equation}
  \label{full_ham}
  H(p,u) = T(p) + U(u), 
\end{equation}
where
\begin{equation}
  \label{eq:15}
  T(p) =  \int_0^1 \tfrac{1}{2} |p(x)|^2 dx, \qquad 
  U(u) =  \int_0^1 \tfrac12\delta^2 |u_x(x)|^2 +V(u(x)) dx.
\end{equation}
We will consider~\eqref{full_pde0} in a parameter regime in which the energy~\eqref{full_ham} is infinite (in which case we will have to consider weak solutions of~\eqref{full_pde0} -- more on this below). In this regime, equation \eqref{full_pde0} possesses a single natural invariant measure which coincides with the canonical (Gibbs) measure of the associated Langevin system, i.e. the stochastic partial differential equation (SPDE) formally given by
\begin{equation}\label{full_damped}
u_t= p, \qquad p_t = \delta^2 u_{xx} - V'(u)  - \gamma p + \sqrt{2\gamma\beta^{-1}}\,\eta(x,t) ,
\end{equation}
in which $\eta(x,t)$ is a space-time white noise, i.e. the Gaussian process with mean zero and covariance $\EE (\eta(x,t)\eta(y,s)) = \delta(x-y)\delta(t-s)$, $\gamma>0$ is the friction coefficient, and $\beta>0$ is the inverse temperature. The Gibbs measure for~\eqref{full_damped} can be written formally and in an unnormalized way using the notations of McKean and Vaninsky~\cite{MV1999} as
\begin{equation}
\label{full_IM0}
d\mu = \exp\left({-\beta H(p,u)}\right) d^\infty p d^\infty u .
\end{equation}
The precise meaning of this measure will be explained in Sec.~\ref{Sec:equivalence}.  The fact that~\eqref{full_IM0} is the only natural invariant measure for~\eqref{full_pde0} is intriguing: indeed, in finite dimensions, Hamiltonian systems possess two different invariant measures, the microcanonical and the canonical ones, and if their dynamics are ergodic, it can only be so with respect to the former since their evolution is energy-conserving~\cite{Sinai_book}. What therefore happens to the microcanonical measure in the context of the wave equation?  The equivalence of invariant measure for the energy conserving wave equation in~\eqref{full_pde0} and its Langevin counterpart in~\eqref{full_damped} also raises another question: Do other features of their dynamics coincide, in particular those related to metastability?

We will address both these questions here. Regarding the first, we will consider finite dimensional truncations of the wave equation based on the Hamiltonian defined as
\begin{equation}
  \label{discrete_ham}
  H_N(\vec{p},\vec{u}) = \frac{1}{N}\sum_{j=1}^{N} \left(\tfrac12p_j^2 + V(u_j)\right)
  + N\sum_{j=0}^{N}\tfrac12 \delta^2(u_{j+1}-u_j)^2,
\end{equation}
where $u_j = u(x_j,t)$ and $p_j = p(x_j,t)$ for $x_j = j/(N+1)$, we denote $\vec{u} = (u_1,\ldots,u_N)$ and $\vec{p} = (p_1,\ldots,p_N)$, and $u_0$ and $u_{N+1} $ are fixed by some boundary conditions.  These finite dimensional truncations do indeed possess both the microcanonical and the canonical measures (see \eqref{micro} and~\eqref{canonical} below) as invariant measures. However we will show in Proposition~\ref{prop:1} below that \emph{both the microcanonical measure on the energy shell $E_N = N /\beta$ and the canonical measure at inverse temperature $\beta$ associated with the truncated Hamiltonian in~\eqref{discrete_ham} converge to the same canonical measure in~\eqref{full_IM0} as the size $N$ of the truncation tends to infinity.}  This result will be proven in the context of quadratic potentials $V(u)$ but should hold also for more general potentials. Note that the scaling $E_N=N/\beta$ implies that the energy of the discretized solutions converge to infinity as $N\to\infty$.

The inverse temperature $\beta$ in~\eqref{full_IM0} controls the concentration of the invariant measure around the minimizers of the Hamiltonian~\eqref{full_ham}.  This concentration of measure suggests that the system's dynamics can display metastability if $V(u)$ is a double well potential, like e.g. $V(u) = \frac14(1-u^2)^2$, since in that case the Hamiltonian~\eqref{full_ham} has two minimizers: these are $(p,u)=(0,u^-)$ and $(p,u)=(0,u^+)$ where $u^\pm$ are the solutions to
\begin{equation}
  \label{eq:13}
  \delta^2 u_{xx}= V'(u) 
\end{equation}
that minimize the potential energy $U(u)$ defined in~\eqref{eq:15}. For $V(u) = \frac14(1-u^2)^2$ and Neumann boundary condition, these are simply $u^\pm = \pm 1$, and they can be similarly obtained for Dirichlet boundary conditions (see Fig.~\ref{fig:solutions} below). We will show in Proposition~\ref{prop:concentration} that the measure in~\eqref{full_IM0} concentrates around these minimizers $(0,u^\pm)$ as $\beta\to\infty$, therefore suggesting that the solutions to~\eqref{full_pde0} spend long periods of time around one before making a transition to the vicinity of the other.

As a means to quantify this behavior and answer the second question raised above, we will compute the mean frequency, averaged over initial conditions drawn from the invariant measure~\eqref{full_IM0}, at which the solutions to~\eqref{full_pde0} make transitions between two regions that each contain $(0,u^-)$ and $(0,u^+)$ and partition its phase-space. Specifically, we will take as the boundary between these regions the affine set of codimension 1 defined as
\begin{equation}
  \label{eq:14}
  \mathcal{S} = \left\{ (p,u): \textstyle \int_0^1 (u-u^s) \phi^s  dx = 0\right\}
\end{equation}
where $u^s(x)$ is the saddle configuration between $u^-$ and $u^+$, i.e. the solution to~\eqref{eq:13} with Morse index 1 that minimizes $U(u)$, and $\phi^s(x)$ is the unique eigenfunction with negative eigenvalue of the Hessian operator at $u^s$, that is
\begin{equation}
  \label{opAs}
  \mathcal{A}^s = -\delta^2\frac{d^2}{dx^2}  + V''(u^s(x)).
\end{equation}
Denoting by $N_T(p_0,u_0)$ the total number of times the solution to~\eqref{full_pde0} for the initial condition $(p(0),u(0))=(p_0,u_0)$ crosses the dividing surface $\mathcal{S}$ during the time interval $[0,T]$, $T>0$, we will compute
\begin{equation}
  \label{eq:12}
 \nu_{\mathcal{S}} = \lim_{T\to\infty} \frac1{2T} \int N_T(p_0,u_0) d\mu(p_0,u_0) \equiv \lim_{T\to\infty} \frac{\EE_\mu N_T}{2T} .
\end{equation}
Transition state theory (TST)~\cite{Eyring1935,Horiuti1938} permits the exact calculation of this quantity: here we will generalize the derivation in~\cite{TalEVE2006} to the wave equation with random initial data to obtain an exact expression for $\nu_{\mathcal{S}}$ that we will then estimate asymptotically in the limit $\beta\to\infty$. This will allow us to obtain asymptotic expressions for the average residency times $\tau_-$ and $\tau_+$ on either side of $\mathcal{S}$, defined as
\begin{equation}
  \label{eq:16}
  \tau_{\pm} = \frac{\mu(\mathcal{B}_\pm)}{\nu_{\mathcal{S}}}
\end{equation}
where 
\begin{equation}\begin{aligned}
  \label{eq:17}
  \mathcal{B}_- &= \left\{ (p,u): \textstyle \int_0^1 (u-u^s) \phi^s  dx < 0\right\}, \quad \mathcal{B}_+ &= \left\{ (p,u): \textstyle \int_0^1 (u-u^s) \phi^s  dx > 0\right\} .
\end{aligned}\end{equation}
The asymptotic expressions for $\tau_\pm$ we will arrive at are
\begin{equation}
  \label{tau}
  \tau_\pm \sim 2 \pi \Lambda_\pm e^{\beta \mathit{\Delta} E_\pm}\qquad 
  \text{as} \ \beta \to\infty .
\end{equation}
Here, $\mathit{\Delta} E_\pm = U(u^s)-U(u^\pm)$ denote the energy barriers between the minima and the saddle point  and 
\begin{equation}
  \label{Lambda}
  \Lambda_\pm =  \frac{1}{\sqrt{ \lambda_1^\pm}} 
  \prod_{j=2}^\infty  \sqrt{\frac{\lambda_j^s}{ \lambda_j^\pm}}
\end{equation}
where $\lambda_j^s$ are the eigenvalues of 
the operator~\eqref{opAs} and $\lambda_j^\pm$ those of
\begin{equation}
\label{opA}
\mathcal{A}^\pm = -\delta^2\frac{d^2}{dx^2}  + V''(u^\pm(x)).
\end{equation}
In~\eqref{Lambda} the eigenvalues are enumerated such that $\lambda_j<\lambda_k$ whenever $j<k$. Note that the negative eigenvalue $\lambda_1^s<0$ corresponding to the unstable direction of the saddle solution (that is, associated with the eigenvector $\phi^s$ defined above) is not included in the infinite product. Note also that the infinite product in \eqref{Lambda} converges in one spatial-dimension, $D=1$, but it does not for $D\ge 2$. This issue is consistent with the behavior of the SPDE in~\eqref{full_damped} which is well-posed for $D=1$, but not for $D\ge2$~\cite{Hairer2014,RNT2012}.  Interestingly the prediction in~\eqref{tau} for the residency times of the solutions to~\eqref{full_pde0} coincide with the zero-damping limit of the residency times of the solutions to the Langevin stochastic partial differential equation in~\eqref{full_damped}, since the latter are given by
\begin{equation}
  \label{tauSPDE}
  \tau^{\gamma}_\pm \sim  \pi \frac{\gamma + \sqrt{\gamma^2 
      + 4|\lambda_1^s|}}{\sqrt{|\lambda_1^s|}} \Lambda_\pm e^{\beta \mathit{\Delta} E_\pm}
  \qquad \text{as \ $\beta \to\infty$ \ with \ $\gamma>0$ fixed.}
\end{equation}

An interesting additional question is whether these results also apply for almost all initial conditions drawn from the invariant measure. This trajectory-wise comparison is nontrivial as it requires an assumption of ergodicity.  Below, we will investigate this assumption numerically. These computations suggest that the system's dynamics is indeed ergodic and rapidly mixing if the parameter $\delta$ is small enough: in this case, the residency times of every trajectory closely match that predicted by TST, and successive visits to $\mathcal{B}_+$ and $\mathcal{B}_-$ are roughly uncorrelated statistically. This means that the dynamics of~\eqref{full_pde0} can be effectively coarse-grained to that of a two-state Markov jump process with rates that can be deduced from the mean TST frequency. It also means that in this regime the dynamics of~\eqref{full_pde0} is not fundamentally different from that of its stochastic counterpart in~\eqref{full_damped}, and in the latter we expect that the limits $\beta\to\infty$ and $\gamma\to0$ commute. If the parameter $\delta$ is too large, however, either ergodicity breaks down, or the mixing time increases and becomes much larger than the inverse of the mean TST frequency. As a result, successive visits in $\mathcal{B}_+$ and $\mathcal{B}_-$ become strongly correlated, and the coarse-graining of~\eqref{full_pde0} to a two-state Markov jump process fails. In this regime, the dynamics of~\eqref{full_pde0} is different from that of~\eqref{full_damped}, and in this second equation, the limits as $\beta\to\infty$ and $\gamma\to0$ no longer commute.

The remainder of this paper is organized as follows: In Sec.~\ref{sec:setup} we begin by setting up the wave equation along with its finite-dimensional truncation that we will study. In Sec.~\ref{Sec:equivalence} we discuss the equivalence of the microcanonical and canonical distributions in infinite dimensions, with detailed calculations of the convergence of characteristic functions of the solutions in appendix~\ref{app:micro_IM} and of the concentration of measure in appendix~\ref{app:quad}.  In Sec.~\ref{sec:TST} we then revisit TST and derive the mean frequency of transition between two specific regions in the phase-space of solutions of the wave equation: these calculations are done in main text for initial conditions drawn from the canonical distribution and in appendix~\ref{app:micro} for initial conditions drawn from the microcanonical distribution. In Sec.~\ref{sec:numerical}, we confirm these results via numerical calculations and also investigate the ergodicity of the dynamics in different parameter regimes.  Finally, in Sec.~\ref{sec:conclusions}, we conclude the paper by discussing the implications of the Hamiltonian dynamics on the low-damping regime of the Langevin system.

\section{Preliminaries: Discretized Equation, Weak Solutions, and Critical Points\label{sec:setup}}

To understand the behavior of the solution to~\eqref{full_pde0}, we will consider the finite discretization of the interval $x\in[0,1]$, into $N+1$ points, defining $x_j = j/(N+1)$, $u_j = u(x_j,t)$ and $p_j = u_t(x_j,t)$, and work with the discrete Hamiltonian in~\eqref{discrete_ham}.  Hamilton's equations associated with~\eqref{discrete_ham} are
\begin{equation}\label{discrete_sys}
\dot{u}_j = p_j, \qquad 
\dot{p}_j = \delta^2 N^2 ( u_{j-1} -2u_j + u_{j+1}) - V'(u_j)
\end{equation}
for $j=1$ to $N$, where the dot indicates differentiation with respect to time.  The Hamiltonian also has $2N$ degrees of freedom as the endpoints $u_0$ and $u_{N+1}$ are fixed by the boundary conditions.  For example, $u_0=u_{N+1}=0$ for Dirichlet boundary conditions, or $u_0=u_1$ and $u_{N+1}=u_N$ for Neumann boundary conditions.  

The solutions to~\eqref{discrete_sys} converge to those of the wave equation in~\eqref{full_pde0} as $N\to\infty$. However, since we will consider a parameter regime in which the energy of the initial conditions to \eqref{discrete_sys} scale as $N/\beta$, the discretized solutions to~\eqref{discrete_sys} will converge to infinite energy solutions of \eqref{full_pde0}, and we cannot expect these solutions to be differentiable in space and time. Therefore, we will have to consider weak solution of~\eqref{full_pde0} (or rather \eqref{full_pde00}) that satisfy
\begin{equation}
  \label{full_pde_weak}
  \int_0^\infty \int_0^1 \left( \phi_{tt} - \delta^2 \phi_{xx} \right) u + V'(u) \phi \;dx dt 
  = \int_0^1 \phi(x,0)u_t(x,0) - \phi_t(x,0)u(x,0) dx
\end{equation} 
for all smooth test functions $\phi(x,t)\in C^\infty$ with compact support in time, the same boundary conditions as $u(x,t)$, and with $\phi_t(x,0)=0$.

As mentioned in the introduction, the critical points of the Hamiltonian in~\eqref{full_ham} will play an important role in what follows. Because of the separable structure of the Hamiltonian, these critical points are always of the form $(p,u)=(0,u)$ with $u$ a solution to~\eqref{eq:13}. The specific shape of these solutions depends both on the form of $V(u)$ and the choice of boundary conditions. To fix ideas, consider e.g.  \begin{equation}
  \label{double_well}
  V(u) = \frac14(1-u^2)^2
\end{equation}
and Neumann boundary conditions.  In this case the solutions of~\eqref{eq:13} that minimize $U(u)$ in~\eqref{eq:15} (i.e. such that all the eigenvalues of the operator \eqref{opA} are positive) are  the constant functions $u^+(x)=1$ and $u^-(x)=-1$.

\begin{figure}[t]
   \centering
   \includegraphics[scale=0.38]{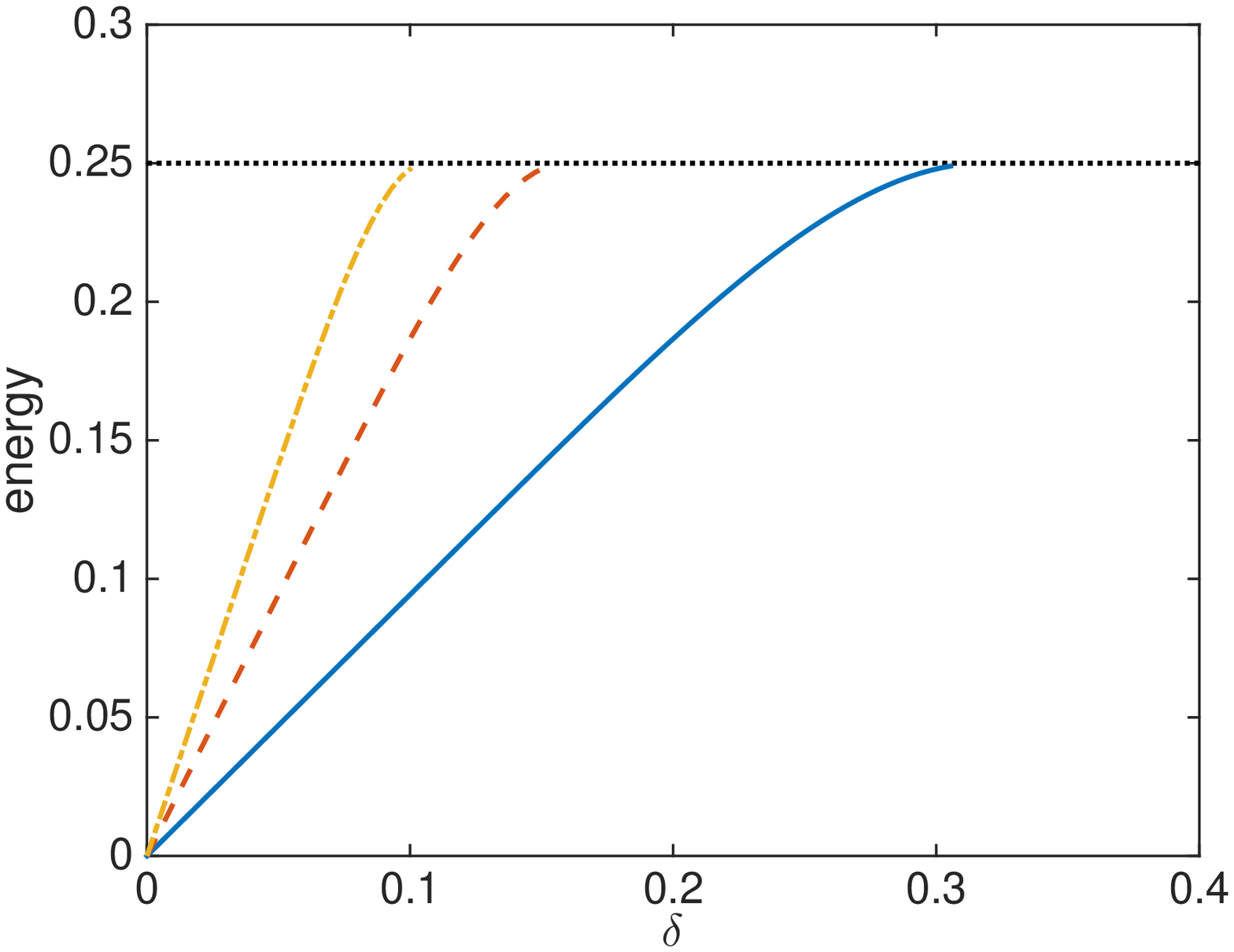} 
    \includegraphics[scale=0.38]{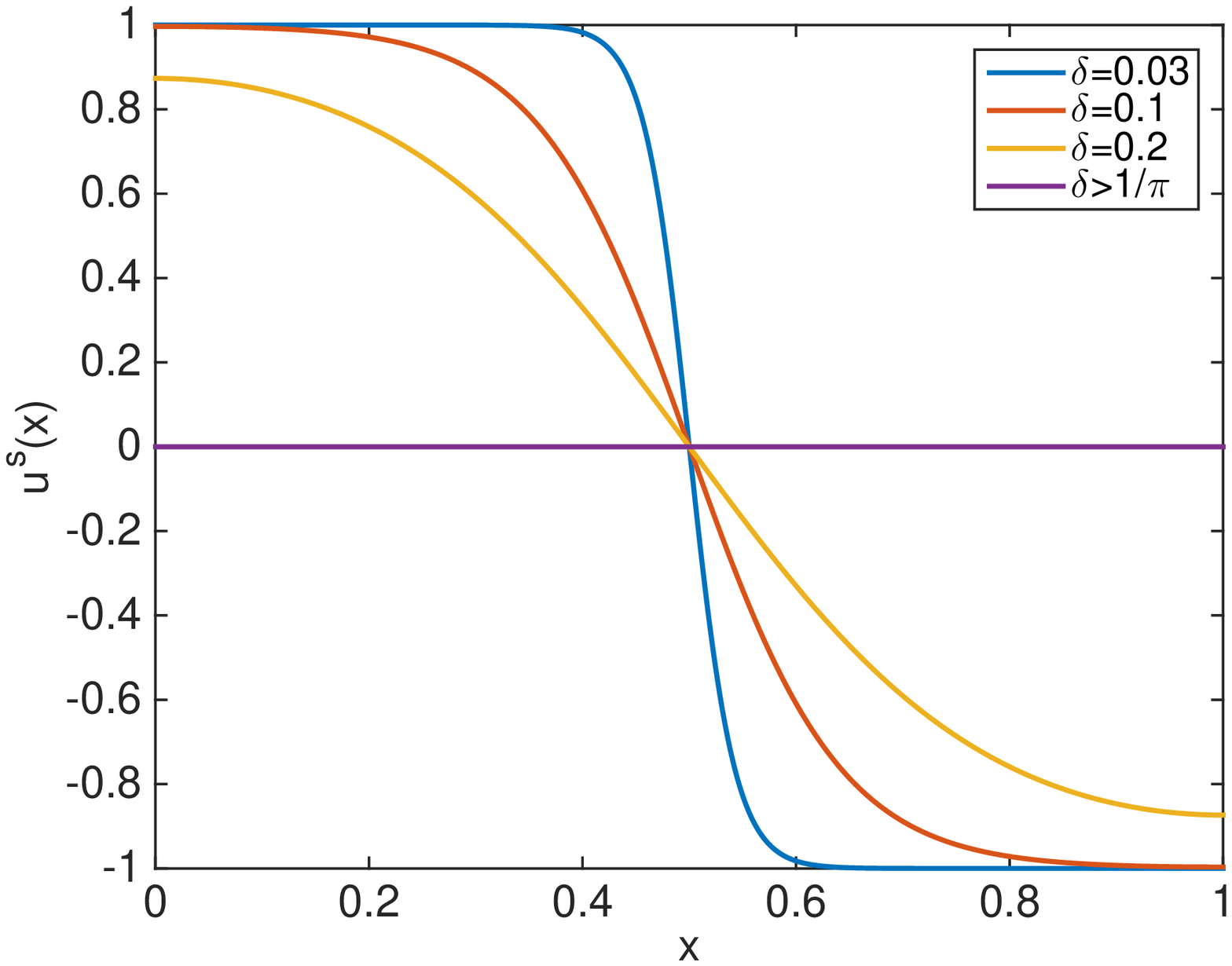}
    \caption{Critical points of $U(u)$ for the potential given in~\eqref{double_well} and Neumann boundary conditions. Left panel: Energies of the first three non-trivial saddle points that emerge as $\delta$ is decreased.  The bifurcations happen at $1/n\pi$ for $n=1,2,3,\dots$ corresponding to the $u(x)=0$ solution with energy $1/4$ gaining negative eigenvalues. Right panel: The index 1 saddle point $u^s(x)$ (one negative eigenvalue) for various values of $\delta$.  For $\delta<1/\pi$ these points correspond to points along the blue solid line in the left panel.  Also for $\delta<1/\pi$, due to symmetry, there is a second index 1 saddle point $u^s(-x)$.}
   \label{fig:bifurcation}
\end{figure}

Besides the minimizers of $U(u)$, we will also need the saddle point with Morse index 1, i.e. the function $u^s(x)$ that is a solution to~\eqref{eq:13} such that the operator \eqref{opAs} has only one negative eigenvalue.  Unlike the minimizers, $u^s(x)$ as well as its energy $E^s$ depend on the value of $\delta$ (see Fig.~\ref{fig:bifurcation} (right)).  For instance, $u(x)=0$ is always a critical point of the Hamiltonian, but it is only the index 1 saddle point when $\delta>1/\pi$.  This can be seen from the eigenvalues of the operator \eqref{opAs} evaluated at $u^s=0$, namely $\lambda = -1+\delta^2n^2 \pi^2$ for $n=0,1,2,3\dots$.  The $n=0$ eigenvalue is always $-1$ regardless of the value of $\delta$.  It is the only negative eigenvalue when $ \delta > 1/\pi$, and $u(x)=0$ is therefore the index 1 saddle point.  The $u(x)=0$ solution becomes an index 2 saddle point when $\delta < 1/\pi$, and an index three saddle when $\delta < 1/2\pi$, and so forth as it gains more negative eigenvalues.  These bifurcation points correspond to the appearance of lower energy, lower index, saddle points.  The energy of these saddle points for the first three bifurcations are shown in Fig.~\ref{fig:bifurcation} (left).

A similar analysis can be made if the potential $V(u)$ is modified (e.g.~if we make it asymmetric) or if the boundary conditions are changed
 (e.g.~if we use Dirichlet boundary conditions, $u(0)=u(1)=0$).

%
%
%
%

\section{Existence and Equivalence of Invariant Measures\label{Sec:equivalence}}

The finite dimensional Hamiltonian system associated with~\eqref{discrete_ham} has two invariant measures: the microcanonical measure supported on the energy shell $H_N(\vec{p},\vec{u}) = E_N$, 
\begin{equation}
  \label{micro}
  dm_N = c_N^{-1} \delta(E_N-H_N(\vec{p},\vec{u})) d\vec{p} d\vec{u} 
  \equiv c_N^{-1} \frac{d\sigma(\vec{p},\vec{u})}{|\nabla H_N(\vec{p},\vec{u})|}
\end{equation}
where $d\sigma(\vec{p},\vec{u})$ denotes the surface element on $H_N(\vec{p},\vec{u})= E_N$ and
\begin{equation}
c_N = \int_{\mathbb{R}^{2N}} \delta(E_N-H_N(\vec{p},\vec{u})) d\vec{p} d\vec{u}\equiv
\int_{H_N(\vec{p},\vec{u})= E_N}\frac{d\sigma(\vec{p},\vec{u})}{|\nabla H_N(\vec{p},\vec{u})|};
\label{eq:2}
\end{equation}
and the canonical measure
\begin{equation}\label{canonical}
dM_N = C_N^{-1}e^{-\beta H_N(\vec{p},\vec{u})}d\vec{p}d\vec{u}
\end{equation}
where
\begin{equation}
  C_N=  \int_{\mathbb{R}^{2N}} e^{-\beta H_N(\vec{p},\vec{u})}d\vec{p}d\vec{u}.
  \label{eq:8}
\end{equation}

We claim that if we set $E_N = N/\beta$ and let $N\to\infty$ with $\beta$ fixed, the microcanonical and canonical measures above converge to the same limit. Here we will establish this result in the context of quadratic potential functions $V(u)$. Specifically, the following proposition is proven in appendix~\ref{app:micro_IM}:
\begin{proposition} \label{prop:1}
 Pick two smooth deterministic functions $s(x)$ and $t(x)$, denote $\vec{s}=(s_1, \ldots,s_N)$ and $\vec{t}=(t_1, \ldots,t_N)$ by defining $s_j=s(x_j)$ and $t_j=t(x_j)$ for $j=1\dots N$, and consider the characteristic functions of the microcanonical and the canonical measures defined respectively as
$$\begin{aligned}
\phi_N^m(\vec s,\vec t) &= \int_{\mathbb{R}^{2N}} 
\exp\left( \frac{i}{N}\vec{s}\cdot\vec p +  \frac{i}{N}\vec{t}\cdot\vec u \right) 
dm_N(\vec{p},\vec{u})\\
\phi_N^M(\vec s,\vec t) &= \int_{\mathbb{R}^{2N}} 
\exp\left( \frac{i}{N}\vec{s}\cdot\vec p +  \frac{i}{N}\vec{t}\cdot\vec u \right) 
dM_N(\vec{p},\vec{u})
\end{aligned}$$
Then if we set $E_N = \beta/N$ and $V(u) = \frac12\alpha u^2$, $\alpha>0$, both $\phi_N^m$ and $\phi_N^M$ have the same limit as $N\to\infty$:
\begin{equation}
  \begin{aligned}
    &\lim_{N\to\infty} \phi_N^m(\vec s,\vec t)  = 
    \lim_{N\to\infty} \phi_N^M(\vec s,\vec t) \\
    &\quad  =\phi(s,t) \equiv \exp\left( -\frac{1}{2\beta}\int_0^1  s^2(x) dx+   
      \int_0^1 \int_0^1  t(x) C(x,y) t(y)  \, dx dy \right),
  \end{aligned}
  \label{eq:18}
\end{equation}
where the covariance function $C(x,y)$ solves
\begin{equation}
  \left( -\delta^2\frac{\partial^2}{\partial x^2} + \alpha \right) C(x,y) = \delta(x-y)
  \label{eq:19}
\end{equation}
with the same boundary conditions as~\eqref{full_pde0}.
\end{proposition}

The functional $\phi(s,t)$ in~\eqref{eq:18} is the characteristic functional of a Gaussian measure, and it is consistent with the stochastic processes $p(x)$ and $u(x)$ being statistically independent and such that: $p(x)$ is a spatial white-noise, scaled by the parameter $\beta$, that is a mean-zero Gaussian process with covariance function
 \begin{equation}
   \mathbb{E}[p(x)p(y)] = \beta^{-1} \delta(x-y);
   \label{eq:9}
\end{equation}
and $u(x)$ is a  mean-zero Gaussian process with covariance function
\begin{equation}
   \mathbb{E}[u(x)u(y)] = C(x,y).
   \label{eq:9b}
\end{equation}
where $C(x,y)$ solves~\eqref{eq:19}. Note that the Gaussian measure whose characteristic function is~\eqref{eq:18} can be written formally as~\eqref{full_IM0} with $V(u) = \frac12\alpha u^2$. 

\begin{figure}[t]
   \centering
   \includegraphics[scale=0.38]{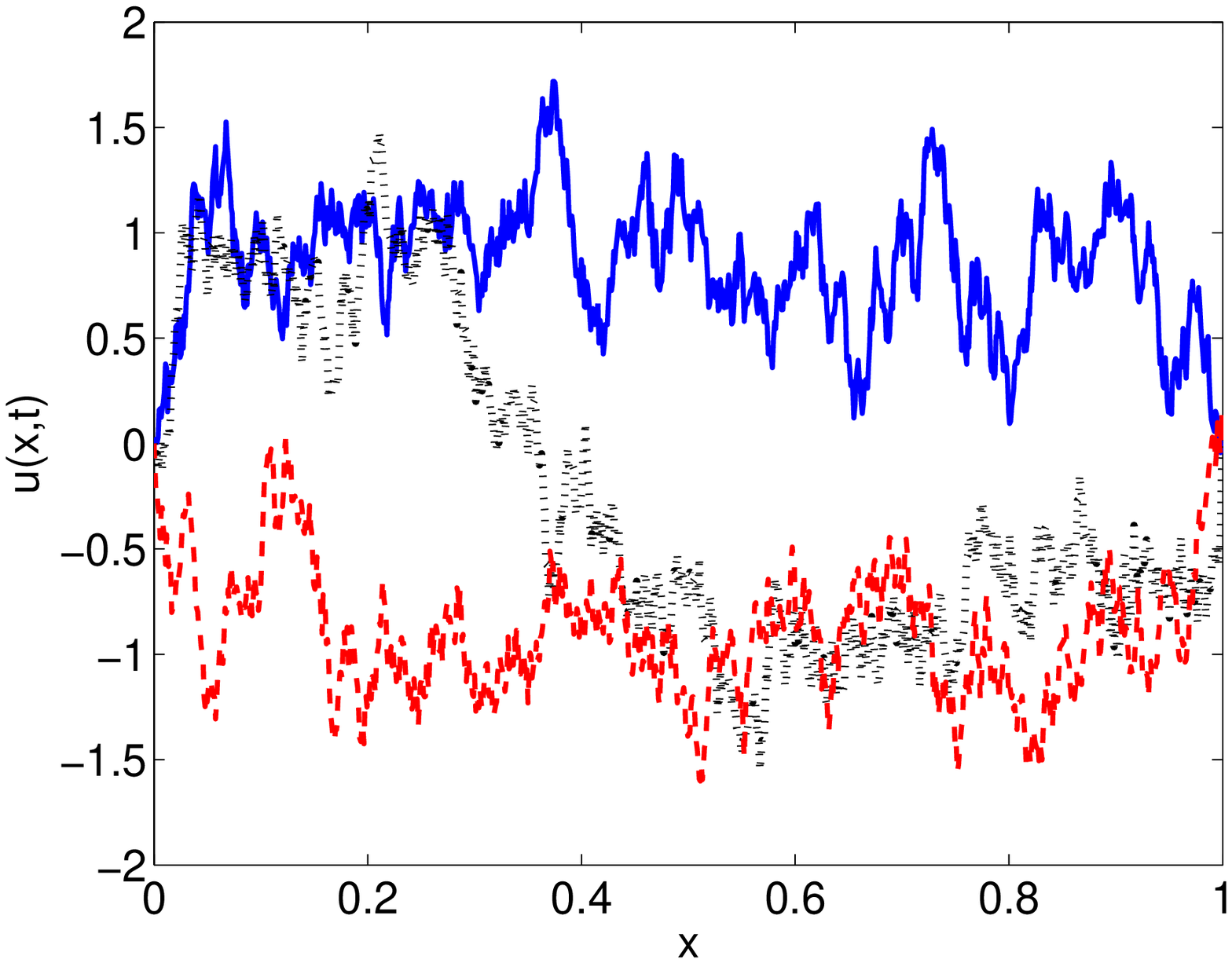} 
   \includegraphics[scale=0.38]{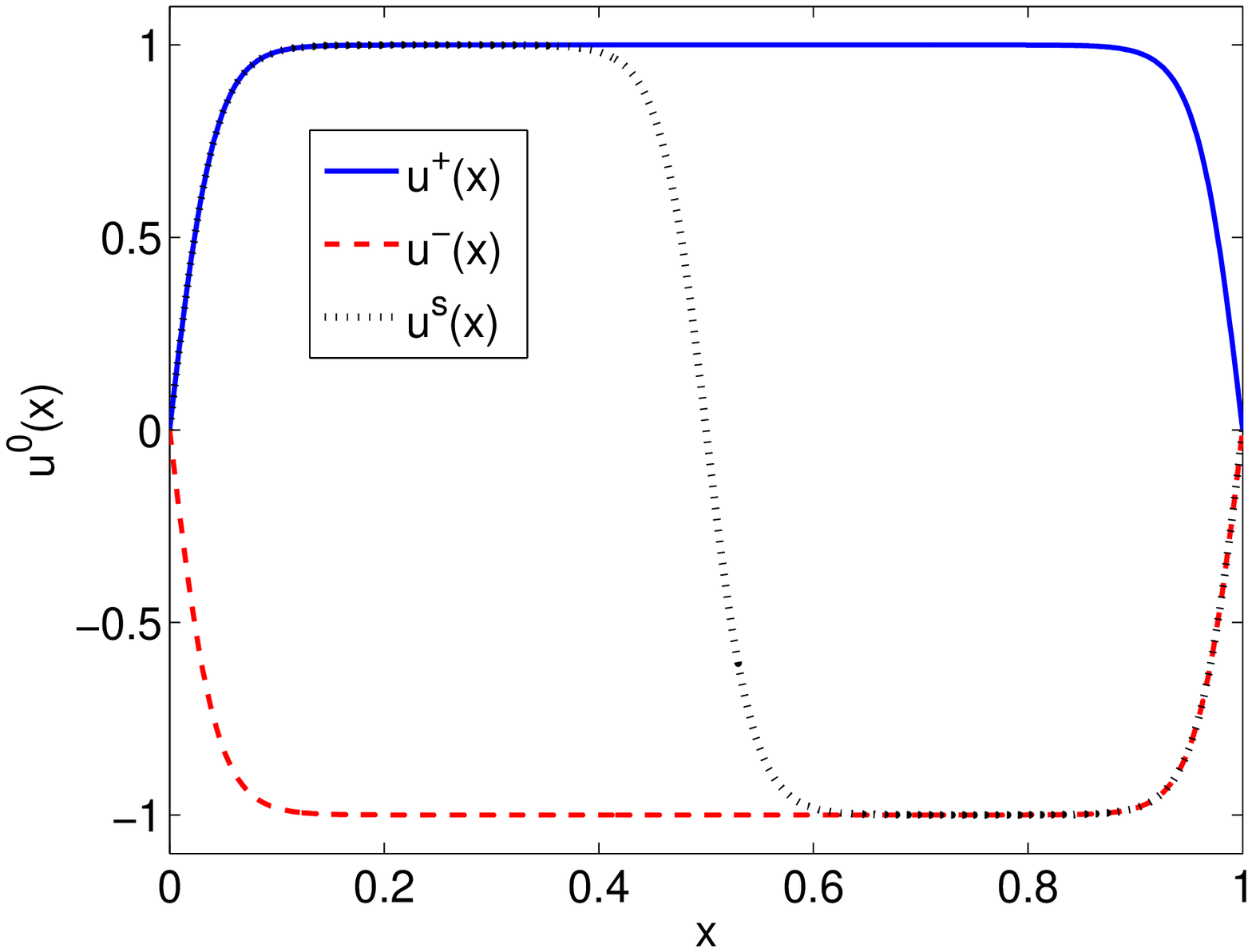}
   \caption{ (Top) Shown are three sample solutions to \eqref{discrete_sys}, the discretized version of \eqref{full_pde0}, with Dirichlet boundary conditions, at various instances in time, with potential function \eqref{double_well}, $N=1024$ discretization points, $\delta = 0.03$, and initial conditions with energy $N/\beta$ and $\beta\Delta E = 4.6$.  The blue solid line is a solution near the minimizer $u^+(x)$, the red dashed line is a solution near the minimizer $u^-(x)$ and the black dotted line is in the process of transitioning from $u^+(x)$ to $u^-(x)$.  (Bottom) These critical points along with one of the saddles, $u^s(x)$, separating them.}
   \label{fig:solutions}
\end{figure}

Proposition~\ref{prop:1} is stated for a quadratic $V(u)$, but we believe that a similar statement holds for a general $V(u)$ that is bounded below and such that $V(u)\ge C|u|^2$ for some $C>0$ as $|u|\to\infty$. Indeed, it seems reasonable to think that the canonical measure is becoming concentrated near the energy surface~$E_N=N/\beta$ in the limit of $N\to\infty$, much in the same way as the volume of a ball becomes concentrated in a thin shell near its surface as the space dimension increases~\cite{Flaschka}. This observation plays a key role in the proof of Proposition~\ref{prop:1} and it suggests that both~\eqref{micro} and~\eqref{canonical} converge to~\eqref{full_IM0} as $N\to\infty$ even for a $V(u)$ that is non-quadratic. In this case, the way to interpret~\eqref{full_IM0} is as follows: the processes $p(x)$ and $u(x)$ are again statistically independent, with $p(x)$ being a spatial white-noise with the covariance in~\eqref{eq:9}, but $u(x)$ is no longer Gaussian.  For example, if we consider Dirichlet boundary conditions, let us introduce the scaled Brownian bridge process $B(x)$, that is, the mean-zero Gaussian process with covariance function \begin{equation}
  \mathbb{E}[B(x)B(y)] = (\beta\delta^2)^{-1}( \min(x,y) - xy).
  \label{eq:10}
\end{equation}
Then $u(x)$ is the the process whose Radon-Nykodim derivative (i.e. density) with respect to $B(x)$ is given by
\begin{equation}
  \exp\left( -\beta\int_0^1 V(B(x)) dx \right).\label{eq:11}
\end{equation}
The parameter $\beta$ appears in this term as well, indicating that $\beta$ and $\delta$ have different influences on the function $u(x)$.  We discuss in Sec.~\ref{sec:quadratic} the effect of $\beta$ on reducing the variance of the solution, equivalently the concentration of the measure, and in Sec.~\ref{sec:numerical} the effect of $\delta$ on the covariance of $u(x)$ leading to more uniform in space solutions.  Shown in Fig.~\ref{fig:solutions} are three snap shots of the solution to \eqref{full_pde0} at three different moments in time, when initial conditions are chosen approximately from the invariant measure \eqref{full_IM0} with the double well potential \eqref{double_well} and Dirichlet boundary conditions.

McKean and Vaninsky~\cite{MV1999} (see also \cite{Friedlander1985,Zhidkov1994}) verified that \eqref{full_IM0} is an invariant measure of \eqref{full_pde0} with periodic boundary conditions $u(0,t) = u(1,t)$.  This invariance of the Gibbs measure is not unique to the wave equation and has been shown for other Hamiltonian PDEs such as the nonlinear Schr\"odinger equation~\cite{Bourgain:1994fs,LEBOWITZ:1988br,Tzvetkov:2006fx,Tzvetkov:2008kh}, Korteweg-de Vries equation~\cite{Oh:2016gm}, and others~\cite{Oh:2010fc,Deng2015}.

We have added to \eqref{full_IM0} the parameter $\beta$ not present in McKean's original results.  Consistent with the aforementioned scaling of the variance of the stochastic functions $p(x)$ and $u(x)$ by $1/\beta$, the parameter $\beta$, added to the measure \eqref{full_IM0}, is playing the roll of an inverse temperature.  Recall the measure \eqref{full_IM0} is also an invariant measure for the noisy-damped system \eqref{full_damped} regardless of the value of $\gamma>0$, the damping parameter.  In this case, the temperature $\beta^{-1}$ characterizes the size of the fluctuations driving the system, and therefore the variance of the solutions.  In terms of the system with Hamiltonian \eqref{full_ham}, the parameter $\beta^{-1}$ acts like a temperature parameter and also characterizes the spatial variance of the weak solutions.

%
%

\section{Expected Residency Times from Transition State Theory\label{sec:TST}}

As a means to determine if other features of the dynamics of the wave equation \eqref{full_pde0} and its Langevin counterpart in \eqref{full_damped} coincide, we will now compute the mean residency times, $\tau_\pm$, the solutions spend in one of two regions that each contain a minimizer of the Hamiltonian and partition phase-space.  To this end, we revisit transition state theory (TST), developed for finite dimensional systems~\cite{Eyring1935,Horiuti1938,Wigner1938}, that computes the mean frequency $\nu_{\mathcal{S}}$ of crossing a codimension one surface partitioning the phase space into two regions; the mean residency times $\tau_\pm$ are related to $\nu_{\mathcal{S}}$ by means of the definition in~\eqref{eq:16}.  Rather than the dynamics of a single trajectory, here we consider the average crossing frequency of a collection of trajectories with initial distribution given by the invariant measure \eqref{full_IM0}.  We start by considering the finite dimensional system~\eqref{discrete_sys}, with initial conditions selected from the invariant measure~\eqref{canonical}, and then justify the infinite dimensional limit.  The same calculation for the microcanonical invariant measure \eqref{micro} is considered in appendix \ref{app:micro}.

\subsection{TST Revisited in an Averaged Set-up}

For a Hamiltonian system, we follow a similar derivation to those in \cite{EVETal2005,TalEVE2006} but do not assume ergodicity at this point.  We divide the finite dimensional phase space, $\vec{u},\vec{p} \in \mathbb{R}^N$ into two regions, $B_-$ containing $(\vec 0,\vec u^-)$ and $B_+$ containing $(\vec 0,\vec u^+)$, by mean of a dividing surface $S$ that separates the two: how to specify this surface will be discussed later in Sec.~\ref{sec:optimalTST}.  For a single trajectory, the average TST frequency of transitions from region $B_-$ to $B_+$ over a time $T>0$ is given by
\begin{equation} 
  \nu_{S}^T(\vec{u}(0),\vec{p}(0)) = \frac{1}{T} N_T(\vec{u}(0),\vec{p}(0)) 
\end{equation} 
where $N_T(\vec{u}(0),\vec{p}(0))$ counts the number of times the trajectory $\vec u(t)$ solution of~\eqref{discrete_sys} for the initial condition $(\vec{u}(0),\vec{p}(0))$ leaves the region $B_-$ by crossing $S$ within the time interval $[0,T]$.  It is convenient to parametrize~$S$ by the zero level set of some function $q(\vec{u})$, specifically,
\begin{equation} 
  S = \{ (\vec p,\vec u): q(\vec{u}) = 0 \}.  
\end{equation} 
We will explain the specific choice of $q(\vec u)$ leading to the surface in \eqref{eq:14} in the next section.  For convention, take $B_+ = \{ (\vec p,\vec u) : q(\vec{u}) > 0 \}$ and $B_-= \{ (\vec p,\vec u) : q(\vec{u}) < 0 \}$.  Transitions from $B_-$ to $B_+$ are then changes in the sign of $q(\vec u)$ from negative to positive.  This allows the frequency to be written as 
\begin{equation} 
  \nu_S^T(\vec{u}(0),\vec{p}(0)) = \frac{1}{T} \int_0^T 
  \max\left( \frac{d}{dt} \Theta \Big(q(\vec u(t))\Big),0\right) dt 
\end{equation} 
where $\Theta(z)$ is the Heaviside function, $\Theta(z)=0$ if $z<0$, $\Theta(0)=\frac12$, $\Theta(z)=1$ if $z>0$. Using the property that $\Theta'(z) = \delta(z)$, where $\delta(z)$ is the Dirac distribution, using the chain rule to get to $\dot{\vec u}(t) \cdot \nabla q(\vec u(t)) \delta(q(\vec u(t)))$, and setting $\dot{\vec u}(t)=\vec p(t)$, we arrive at
\begin{equation}
  \label{eq:k} \nu_S^T(\vec{u}(0),\vec{p}(0)) = \frac{1}{T} \int_0^T \max\left( \vec p(t) \cdot \nabla 
    q(\vec u(t)) ,0\right) \delta(q(\vec u(t))) dt.  
\end{equation}

We now want to compute the average of $\nu_S^T$ over trajectories with initial conditions chosen with respect to an invariant measure, specifically the canonical  distribution \eqref{canonical}.  By definition, if we denote by $(\vec p(t,\vec p_0,\vec u_0),\vec u(t,\vec p_0,\vec u_0))$ the solution to~\eqref{discrete_sys} for the initial condition $(\vec p(0,\vec p_0,\vec u_0),\vec u(0,\vec p_0,\vec u_0)) = (\vec p_0,\vec u_0)$, for any suitable test function $f(p_0,u_0)$, we have
\begin{equation}
  \label{eq:24}
  \int f (\vec p(t,\vec p_0,\vec u_0),\vec u(t,\vec p_0,\vec u_0)) dM_N(\vec p_0,\vec u_0) 
  = \int f (\vec p,\vec u )dM_N(\vec p,\vec u) .
\end{equation}
As a result the ensemble average of the integrand in~\eqref{eq:k} is invariant in time, and after proper interpretation of integration with respect to the Dirac distribution, we arrive at the following expression for the mean frequency
\begin{equation}
  \label{eq:om2_c}
  \nu_S^c =  C_N^{-1} \int_{S}  \max\left( \vec p \cdot \hat n(\vec u),0\right) 
  e^{-\beta H_N(\vec u, \vec p)} d\tilde{\sigma}(\vec p, \vec u) 
\end{equation}   
where $d\tilde{\sigma}(\vec p, \vec u)$ is the surface element on $S$.  Notice that \eqref{eq:om2_c} is  identical to the expression we would have arrived at if we had assumed ergodicity with respect to the invariant measure, and then taken the long time limit of \eqref{eq:k}, replacing the longtime average with an average over the invariant measure.

\subsection{Optimal TST\label{sec:optimalTST}}

To understand metastability, we are interested in evaluating \eqref{eq:om2_c} for the dividing surface $S$ that best separates the two metastable states containing $\vec u^-$ and $\vec u^+$, minimizing the mean frequency of transition~\cite{EVETal2005} (and thereby maximizing the transition time).  
While we may not be able to determine this surface exactly, a reasonable surface to use is the finite dimensional approximation of the codimension one hyperplane given by \eqref{eq:15}.  The hyperplane in \eqref{eq:15} passes  through the saddle point, perpendicular to the unstable direction of the gradient flow
$
w_t(x,t) = \delta^2 w_{xx}(x,t) - V'(w(x,t)).
$
While this surface may not be the true surface that maximizes the mean transition time, for large $\beta$ we expect this surface to have minimal re-crossings, and therefore produce a close bound on the mean transition time.   

The reasoning for this surface having minimal re-crossings starts with the assumption that most transitions will occur near the saddle point, the lowest energy point on the dividing surface~\cite{AEVE2007}.  Near the saddle point, the linearized Hamiltonian system has one hyperbolic direction corresponding to the one unstable direction of the gradient flow; the remaining directions are oscillatory.  Any surface not perpendicular to the unstable direction of the gradient flow (which bisects the stable and unstable hyperbolic directions) would also not be orthogonal to all the oscillatory directions of the Hamiltonian flow, allowing for many recrossings of the surface in at least one of the oscillatory directions.

The finite dimensional approximation to \eqref{eq:15} is the hyperplane
\begin{equation}
  \label{eq:20}
  S = \{(\vec p,\vec u): (\vec u - \vec u^s)\cdot \vec \phi^{(1)} =0\}
\end{equation}
where $\vec u^s$ and $\vec \phi^{(1)}$ are the discrete equivalents to $u^s$ and $\phi^s$, that is:  $\vec u^s$ is the Morse index 1 critical point of
\begin{equation}
  \label{eq:23}
  U_N(\vec u)  = \frac{1}{N}\sum_{j=1}^{N} V(u_j)
  + N\sum_{j=0}^{N}\tfrac12 \delta^2(u_{j+1}-u_j)^2,
\end{equation}
with minimum energy, and $\vec \phi^{(1)}$ is the unique eigenvector with negative eigenvalue of the Hessian of~\eqref{eq:23} evaluated at $\vec u = \vec u^s$. In other words, these vectors solve
\begin{align}
  \label{eq:21}
    0 & = \delta^2 N^2 ( u^s_{j-1} -2u^s_j + u^s_{j+1}) - V'(u^s_j)\\
  \label{eq:22}
   \lambda^s_k \phi_j^{(k)}& = 
    -\delta^2 N^2 ( \phi^{(k)}_{j-1} -2\phi^{(k)}_j + \phi^{(k)}_{j+1}) +V''(u^s_j)\phi^{(k)}_j, 
\end{align}
when $k=1$ and $\lambda^s_1<0$ for $j=1,\ldots, N$ and appropriate boundary conditions fixing $u^s_0$, $u^s_{N+1}$, $\phi^{(k)}_0$, and $\phi^{(k)}_{N+1}$.  Note the remaining eigenvectors for $k=2\dots N$ solve \eqref{eq:22} when $\lambda^s_k>0$.

\subsection{Quadratic Approximation as $\beta\to\infty$\label{sec:quadratic}}

In order to proceed with computing the TST frequency~\eqref{eq:om2_c}, we must perform an integration over the hyperplane given by \eqref{eq:20}, as well as another  integration to evaluate the normalization constant.  While we cannot evaluate the integrals for an arbitrary potential function we can proceed if it is justified to approximate the potential $V(u)$ in the Hamiltonian locally as a quadratic function.  This approximation is justified if most of the measure we are integrating over is confined within a localized region, within which $V(u)$ is approximately quadratic.  We show this concentration for the quadratic potential and work under the assumption that the total measure outside this region for the true potential $V(u)$ is exponentially small,  and therefore even large errors in the approximation of $V(u)$ do not contribute significantly to the integration for the TST frequency.

For the canonical measure, \eqref{canonical}, intuitively one would expect the measure to concentrate around the points $(\vec 0, \vec u^\pm)$ that minimize the Hamiltonian as $\beta\to\infty$.  It is convenient to write an expansion for $\vec u$ in the rotated and rescaled eigenvector basis as
\begin{equation}
\label{eq:45}
\vec u = \vec u^\pm + \tilde{\vec u} = \vec u^\pm +\sum_{k=1}^N b^\pm_k \vec \psi^{(k)\pm},
\end{equation}
in which the eigenvalues $\lambda^\pm_k$ and the corresponding eigenvectors $\vec \psi^{(k)}$ for each $k=1\dots N$ solve
\begin{equation}
\label{eq:46}
\lambda^\pm_k \psi^{(k)\pm}_j = -\delta^2 N^2( \psi^{(k)\pm}_{j-1} - 2\psi_j^{(k)\pm} + \psi^{(k)\pm}_{j+1}) + V''(u^\pm_j)\psi^{(k)\pm}_j
\end{equation}
for $j=1\dots N$ and appropriate boundary conditions fixing $\psi_0^{(k)\pm}$ and $\psi_{N+1}^{(k)\pm}$.  The eigenvectors are normalized to satisfy $\frac{1}{N}\sum_{j=1}^N (\psi^{(k)\pm}_j)^2=1$,  the discrete version of $L_2$ normalization.  We use the same numbering convention as before, $\lambda^\pm_j < \lambda^\pm_k$ if $j<k$, but this time all the $\lambda^\pm_k >0$ for $k=1\dots N$.
The quadratic expansion of the potential energy \eqref{eq:23} is 
\begin{equation}
\begin{aligned}
U_N(\vec u) & \sim U_N(\vec u^\pm) + \frac{1}{2N} \sum_{k=1}^N \lambda^\pm_k (b^\pm_k)^2  \vec \psi^{(k)\pm} \cdot \vec \psi^{(k)\pm}   \quad \textrm{as} \ \vec u \to \vec u^\pm \\ 
& = E_N^\pm + \frac{1}{2} \sum_{k=1}^N \lambda_k^\pm (b^\pm_k)^2  .
\end{aligned}\end{equation}
Defining $\tilde{p}_j = \dot{b}^\pm_j$, the rotated momenta, the quadratic approximation of the Hamiltonian is
\begin{equation}
\label{eq:47}
H_N \sim E_N^\pm + \sum_{j=1}^N \left( \frac{1}{2N} \tilde{p}_j^2 + \frac{1}{2} \lambda^\pm_j (b_j^\pm)^2 \right)  \quad \textrm{as} \ \vec u \to \vec u^\pm,
\end{equation}
where we have defined $E_N^\pm = U_N(\vec u^\pm)$.  We show in appendix~\ref{app:concentration_c}  the canonical measure  concentrates into a localized region around the point $\vec b^\pm = 0$, even as $N\to\infty$ in the following proposition:
\begin{proposition}
\label{prop:concentration}
Consider the marginal canonical measure with the quadratic Hamiltonian approximation \eqref{eq:47}
and an $N$-dimensional box with side edge lengths $2\delta_j$, 
$$
D_\delta^N = \{ \vec b^\pm : -\delta_j \le b^\pm_j \le \delta_j\;\; \forall \;\;j=1\dots N\},
$$
where $\delta_j$ are such that $ \sum_{j=1}^N \delta_j^2$ converges as $N\to\infty$. Assume that the eigenvalues $\lambda^\pm_j$ solving \eqref{eq:46} satisfy
$$
\lambda^\pm_j \sim Cj^2 \qquad \text{as $j\to\infty$, for some $C>0$}
$$
Then, for every $\hat{\delta}>0$, there exists $\beta_1>0$ and a set of $\{ \delta_j\}_{j=1}^N $ such that for every $\beta>\beta_1$ and every $N>0$
\begin{equation}
\frac{\int_{D_\delta^N} \exp\left( -\frac{\beta}{2}\sum_{j=1}^N \lambda^\pm_j (b^\pm_j)^2\right) d\vec b^\pm }
{\int_{\mathbb{R}^N} \exp\left( -\frac{\beta}{2}\sum_{j=1}^N \lambda^\pm_j (b^\pm_j)^2\right) d\vec b^\pm} > 1- \hat{\delta}
\end{equation}
\end{proposition}

Together with the assumption that the total measure outside the localized region is exponentially small, the quadratic approximation \eqref{eq:47} can be used to approximate the normalization constant in~\eqref{eq:om2_c}.
 We first separate the integration into the two regions, $B_+=\{ (\vec p,\vec u): (\vec u-\vec u^s)\cdot \vec\phi^{(1)} >0\}$ and $B_-=\{ (\vec p,\vec u): (\vec u-\vec u^s)\cdot \vec\phi^{(1)} <0\}$ that partition space, separated by $S$, 
\begin{equation}
\label{eq:43}
 C_N = C_N^+ + C_N^-
\end{equation}
where
\begin{equation}
\label{eq:44b}
C_N^\pm = \int_{B_\pm} e^{-\beta H_N(\vec{p},\vec{u})}d\vec{p}d\vec{u}.
\end{equation}
Now, each region contains only one minimizer of the potential energy.  Within region $B_+$, we expect most of the measure is near $\vec u^+$ while in $B_-$, its in the vicinity of $\vec u^-$.  The result after completing the integration appears below in Proposition~\ref{prop:asymptotic}.

For the integral in~\eqref{eq:om2_c}, the measure is restricted to the surface $S$ defined in \eqref{eq:20}.  On this surface, it is the point $(\vec 0,\vec u^s)$  that minimizes the Hamiltonian.  We expand about this point, and rotate to an eigenvector basis 
\begin{equation}
\label{eq:40}
\vec u = \vec u^s + \tilde{\vec u} = \vec u^s +  \sum_{k=1}^N a_k \vec \phi^{(k)} 
\end{equation}
with $\lambda^s_k$ and $\vec \phi^{(k)}$ for each $k=1\dots N$ solving the system in \eqref{eq:22} with appropriate boundary conditions fixing  $\phi^{(k)}_0$, and $\phi^{(k)}_{N+1}$, and normalized to satisfy $\frac{1}{N}\sum_{j=1}^N (\phi^{(k)}_j)^2=1$.
Our convention is to order the eigenvalues such that $\lambda^s_1$ is the sole negative eigenvalue and $\lambda^s_j< \lambda^s_k$ for $j<k$.
The quadratic approximation of the potential energy \eqref{eq:23} is 
\begin{equation}
\begin{aligned}
U_N(\vec u) &\sim U_N(\vec u^s) + \frac{1}{2N} \sum_{k=1}^N \lambda^s_k a_k^2  \vec \phi^{(k)} \cdot \vec \phi^{(k)} \quad \textrm{as} \ \vec u \to \vec u^s \\
& = E_N^s + \frac{1}{2} \sum_{k=1}^N \lambda_k^s a_k^2 
\end{aligned}\end{equation}
due to the defined normalization of the eigenvectors, and where we have defined $E_N^s=U_N(\vec u^s)$.  Defining $\tilde{p}_j = \dot{a}_j$, the rotated momenta,
the approximation of the Hamiltonian in the vicinity of $\vec u^s$ is
\begin{equation}
\label{eq:42}
H_N \sim E_N^s + \sum_{j=1}^N\left( \frac{1}{2N} \tilde{p}_j^2 + \frac{1}{2} \lambda^s_j a_j^2 \right) \quad \textrm{as} \ \vec u \to \vec u^s.
\end{equation}
In appendix \ref{app:concentration_c} we state Proposition \ref{lemma:c_S}, an analogous Proposition to \ref{prop:concentration}, for the concentration of the measure restricted to the surface $S$.  This allows the use of \eqref{eq:42} to complete the integration in \eqref{eq:om2_c}.

The following proposition states the asymptotic in $\beta$ values of the integrals needed to compute the TST frequency in \eqref{eq:om2_c}.  The justification of these formulas appear in appendix~\ref{app:integrals}.
\begin{proposition}
  \label{prop:asymptotic}  
For the canonical ensemble,  the asymptotic in $\beta\to\infty$, valid for all $N>0$, expansion of the integral in \eqref{eq:om2_c}  is
\begin{equation}
\int_{S}  \max\left( \vec p \cdot \hat n(\vec u),0\right) 
  e^{-\beta H_N(\vec u, \vec p)} d\tilde{\sigma}(\vec p, \vec u) \sim  \frac{2^{N-2} N^{N} (N-2)! S_{2N-3}}  {\beta^{N}\prod_{j=2}^N\sqrt{\lambda^s_j}} e^{-\beta E_N^s} 
\end{equation}
and the expansion of the integrals in \eqref{eq:44b} are
\begin{equation}
C_N^\pm  \sim   \frac{  2^{N-1} N^{N} (N-1)!S _{2N-1} }{ \beta^N \prod_{j=1}^N \sqrt{\lambda^\pm_j}} e^{-\beta E_N^\pm}. 
\end{equation}
In the above, the eigenvalues $\lambda_j^s$, $j=2\dots N$ solve \eqref{eq:22}, the eigenvalues $\lambda_j^\pm$, solve \eqref{eq:46}
and $S_{n}=2\pi^{(n+1)/2}/\Gamma((n+1)/2)$ is the surface area of the $n$-dimensional sphere (embedded in $n+1$ dimensional space) of unit radius.
\end{proposition}
Using the above asymptotic in $\beta$ expansions, the TST frequency from \eqref{eq:om2_c} is
\begin{equation}
\nu_S^c \sim \frac{1}{2\pi } \frac{\prod_{j=2}^N \frac{1}{ \sqrt{\lambda^s_j}} e^{-\beta E_N^s}}{\prod_{j=1}^N\frac{1}{ \sqrt{\lambda^+_j}}e^{-\beta E_N^+} + \prod_{j=1}^N \frac{1}{\sqrt{\lambda^-_j}}e^{-\beta E_N^-}}.
\end{equation}

\subsection{Mean Residency Times}
\label{sec:meanrestimes}

The mean residency times on either side of $S$, given by \eqref{eq:20}, can be computed with the same integrals considered above.  Namely, 
\begin{equation}
\tau_\pm^c =  \frac{C_N^\pm}{C_N \nu_S^c}
\end{equation}
for the canonical ensemble.    With the asymptotic expansions of the integrals in Proposition \ref{prop:asymptotic}, we  arrive at the asymptotic in $\beta$ expressions for the mean residency time
\begin{equation}\label{tau_finite3}
\tau_\pm^c   \sim 2 \pi \frac{\prod_{j=2}^N \sqrt{\lambda_j^s}}{\prod_{j=1}^N \sqrt{\lambda_j^\pm}} e^{ \beta (E_N^s - E_N^\pm)} .
\end{equation}
This expression is valid for all $N$ allowing us to take the limit as $N\to\infty$, arriving at $\tau_\pm$ in \eqref{tau}, provided the  
 infinite product $\Lambda_\pm$ give in \eqref{Lambda} converges.  It is the infinite dimensional limit of the finite dimensional eigenvalue products in    \eqref{tau_finite3}.
 For example, $\Lambda_\pm $ converges when considering eigenvalues related to the Laplacian in one spatial dimension, as we do here (however it does not for two or more spatial dimensions).

\section{Numerical Results\label{sec:numerical}}
 
Using numerical simulations of the finite dimensional truncation of \eqref{full_pde0}, we investigate the implication of our results for the long time dynamics of the nonlinear wave equation. We proceed in two steps. First, we validate our approximations by considering an ensemble of trajectories chosen approximately from the finite dimensional truncation of the invariant measure~\eqref{full_IM0}, and estimating numerically the mean residency time averaged over these initial conditions.  Second, we compute the residency time averaged over a single trajectory in order to investigate the extent to which the dynamics of \eqref{full_pde0} is ergodic and mixing with respect to the invariant measure~\eqref{full_IM0}. We find support for this statement for small values of the exchange coupling $\delta$.  For larger values of $\delta$, however, we observe that either ergodicity fails, or the mixing time grows so much that we cannot observe convergence of the time averages towards their ensemble counterparts averaged over the invariant measure~\eqref{full_IM0}.

We run numerical simulations of the discretized system of $2N$ first order differential equations~\eqref{discrete_sys} with Neumann boundary conditions using the 
velocity-Verlet integrator to conserve the energy on average.  Initial conditions are chosen from an approximation of the invariant measure \eqref{full_IM0}: the finite dimensional canonical measure \eqref{canonical} with the quadratic approximation for the Hamiltonian~\eqref{eq:47}.  Specifically,  $\vec p \in \mathbb{R}^N$ is a discretization of white noise and 
$\vec u \in \mathbb{R}^N$ is a Gaussian random variable with known mean and covariance.  They are computed as
\begin{equation}
\label{IC}
\vec{p} = \sqrt{\frac{N}{ \beta}}\; \vec{h} \; , \quad
\vec{u} = \vec{u}^+ + \frac{\alpha}{\sqrt{\beta}} \sum_{j=1}^N g_j \sqrt{\lambda_j^+}\vec \psi^{(k)+}
\end{equation}
in which $\vec h=(h_1,\dots h_N)$ and $\vec g = (g_1,\dots g_N)$ are vectors of $N$ independent mean zero standard Gaussian random variables,  $\vec u^+=1$ for Neumann boundary conditions, the eigenvalues $\lambda_+$ and eigenvectors $\vec \psi^{(k)+}$ solve \eqref{eq:46},  normalized such that $\vec \psi^{(k)+}\cdot \vec \psi^{(k)+} = N$, and $\alpha$ is chosen to ensure the energy of the initial condition is $N/\beta$ exactly.

For the potential function \eqref{double_well}, the system is symmetric.  We therefore determine the mean transition time by taking the total simulation time and divide it by the number of transitions across the diving hyperplane $S$ described in section \ref{sec:optimalTST}.  Specifically, we compute the eigenvector $\vec \psi^{(1)}$ solving \eqref{eq:22}, the normal vector to the surface $S$.  At each time step, we compute the distance to this dividing surface.  A change in sign of the distance indicates a crossing of the surface.  Notice that for $\delta < 1/\pi$ there are two saddle points, yet they both lie in the just described plane.  

\begin{figure}[t]
   \centering
   \includegraphics[scale=0.38]{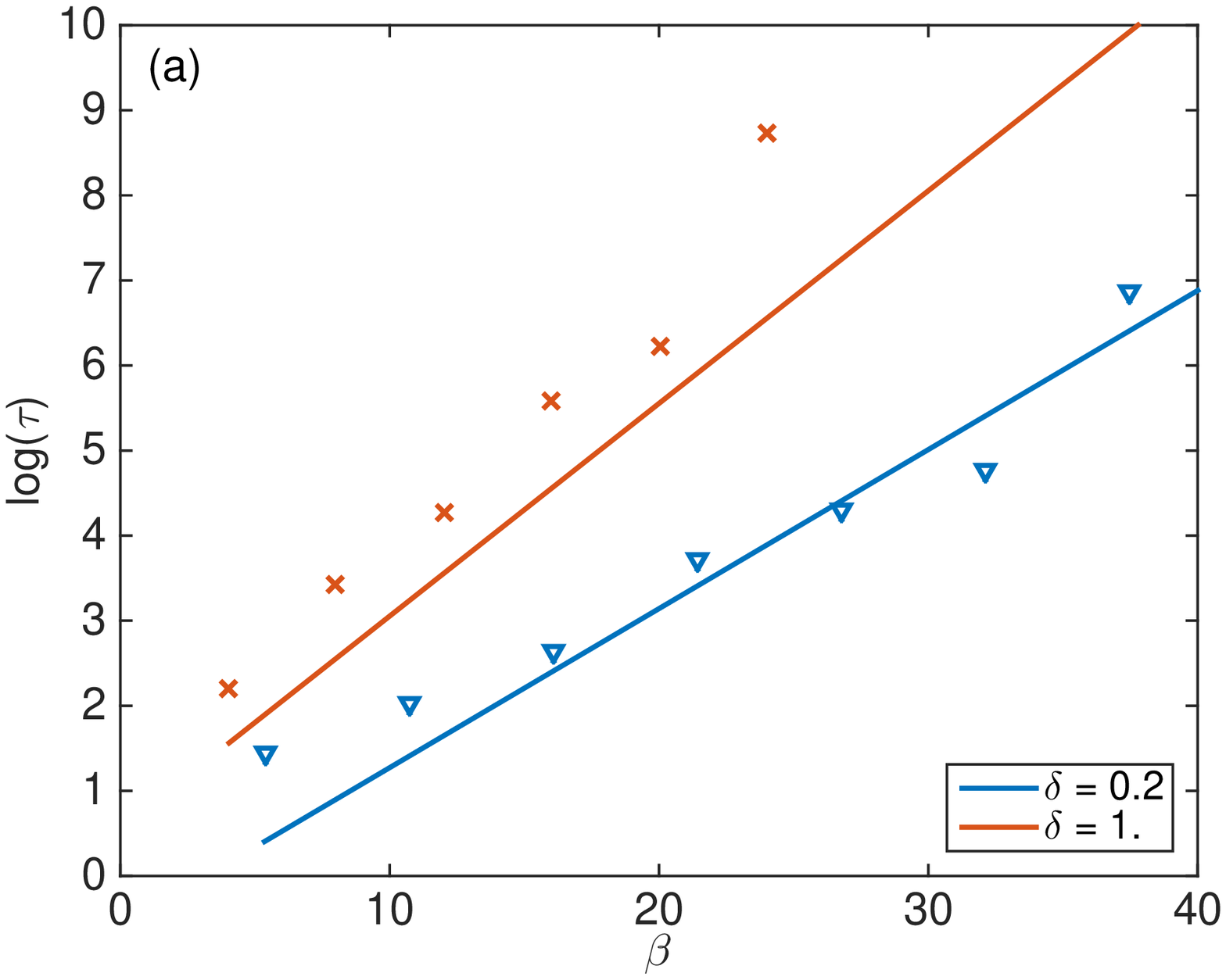}
   \includegraphics[scale=0.38]{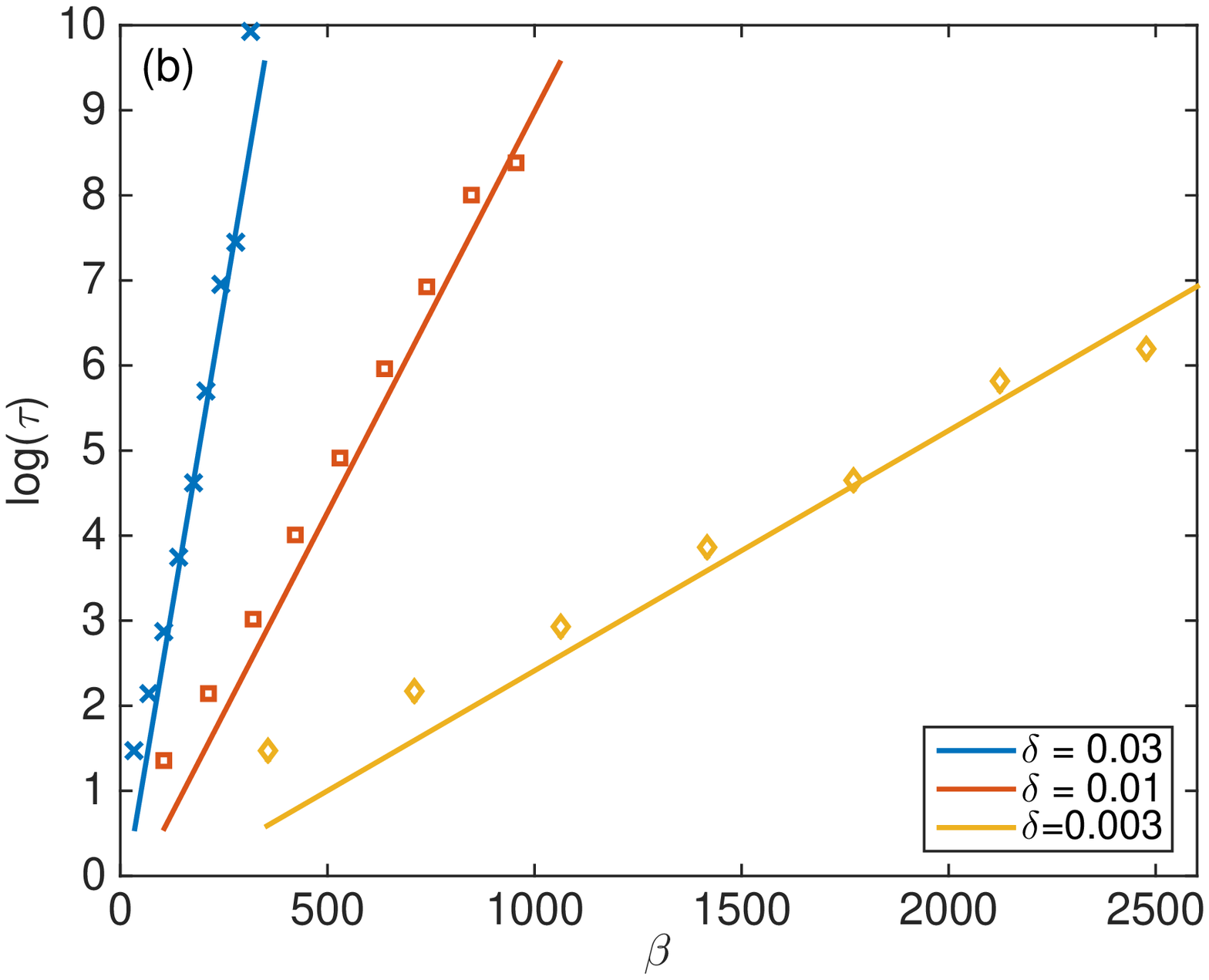} 
      \includegraphics[scale=0.38]{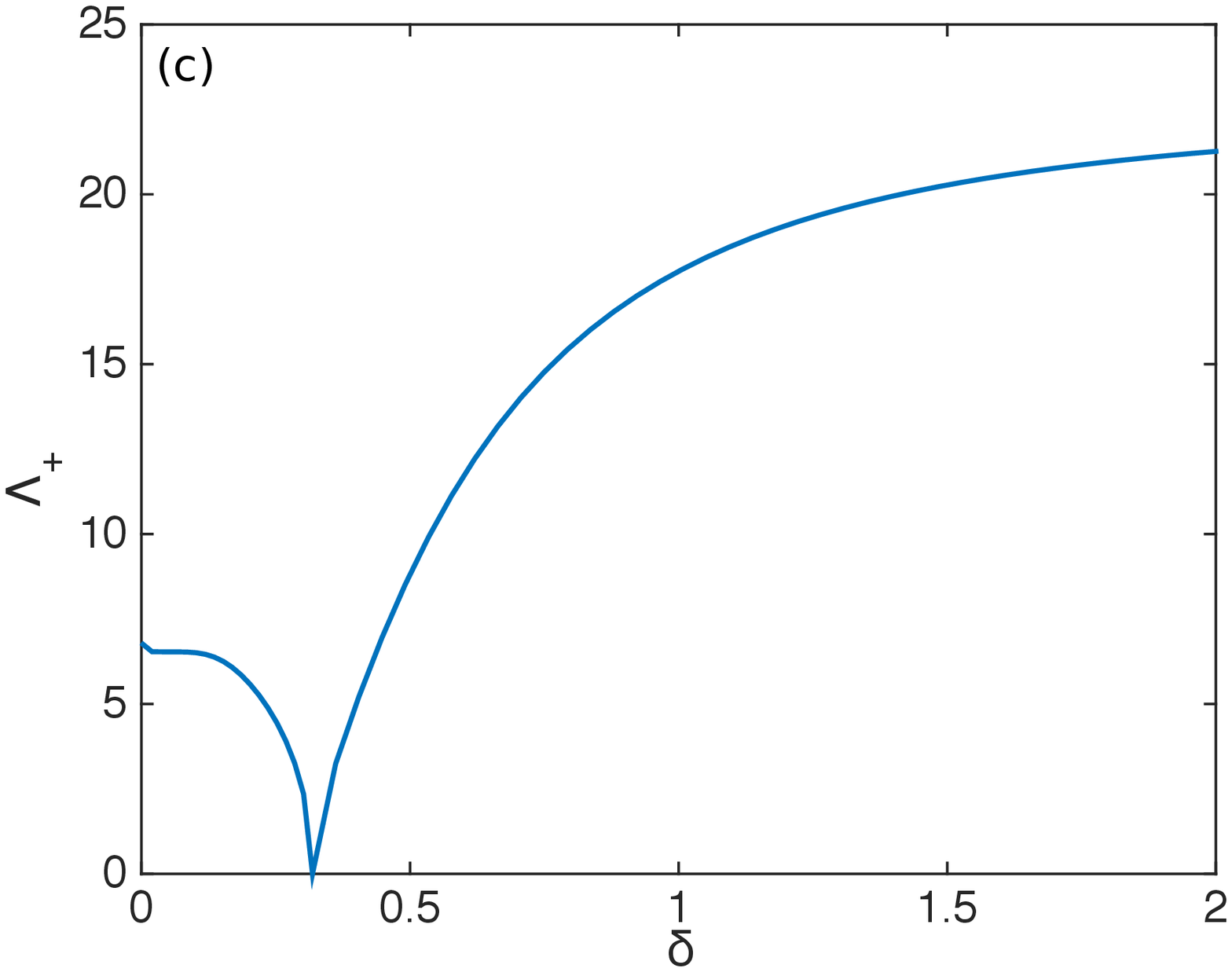}
   \caption{(a) Log of the mean transition time as a function of $\beta$ (symbols) averaged over an ensemble of 1000 trajectories chosen according to \eqref{IC}, evolved by numerically simulating
    \eqref{discrete_sys}, with potential function \eqref{double_well} and $N=128$ and (lines) the theoretical value of $\tau_+$ given by \eqref{tau}. (b) Log of the mean transition time (symbols) computed from a long time average of a single trajectory chosen according to \eqref{IC} and (lines) the theoretical value of $\tau_+$ given by \eqref{tau}.  
   (c)   The prefactor $\Lambda_+$ in \eqref{Lambda} computed numerically with $N=1024$ points as a function of $\delta$. As the bifurcation point $\delta = 1/\pi$ is approached, the prefactor goes to zero because one of the eigenvalues of the saddle point is about to cross zero.  
   }
   \label{fig:tau}
\end{figure}

In Fig.~\ref{fig:tau} (a) we show excellent agreement between the mean transition time found from the numerical simulations at finite $N=128$ averaged over an ensemble of 1000 trajectories, and the infinite dimensional theory, \eqref{tau}.  We choose two values of $\delta$, one smaller than the bifurcation point $\delta=1/\pi$ in which $u^s$ is non-uniform and one larger than this point in which case $u^s=0$.  This agreement supports the validity of the quadratic approximation used in obtaining \eqref{tau}.

Next we investigate the ability of the TST time to describe the mean transition time of a single trajectory.  As we pointed out above in the derivation of $\tau_\pm$, we would arrive at the same expression had we assumed the flow was ergodic with respect to the invariant measure.  Recall in finite dimensions, no Hamiltonian systems are ergodic with respect to the invariant measure given by the canonical ensemble while many are conjectured to be ergodic with respect to the invariant measure given by the micro-canonical ensemble~\cite{Sinai_book}.  In infinite dimensions, McKean and Vaninsky conjectured that the flow of \eqref{full_pde0} is ergodic (a.k.a. metrically transitive) with respect to the measure \eqref{full_IM0} except for completely integrable systems with an infinite number of conserved quantities like $V'(u)=u, \sin u$ or $ \textrm{sinh } u$~\cite{MV1999}.  Although the measure \eqref{full_IM0} looks like the canonical distribution, we argued its equivalence to the microcanonical measure, and therefore do not find an assumption of ergodicity surprising.

Indeed, at small values of $\delta$, we find numerical evidence for ergodicity, and agreement between the mean transition time found from a long-time average over a single trajectory and the time \eqref{tau} as shown in Fig.~\ref{fig:tau}(b).   In this semi-log plot, the slope of the lines are $\mathit{\Delta} E_+ $ (for Neumann boundary conditions this is solely $E^s$, the value of which is show in Fig.~\ref{fig:bifurcation} as a function of $\delta$), while the y-intercept is given by the prefactor, $2\pi \Lambda_+$.  The eigenvalue produce $\Lambda_+$ is shown in Fig.~\ref{fig:tau}(c) as a function of $\delta$, computed with a truncation of $N=1024$.

\begin{figure}[t]
   \centering
         \includegraphics[scale=.38]{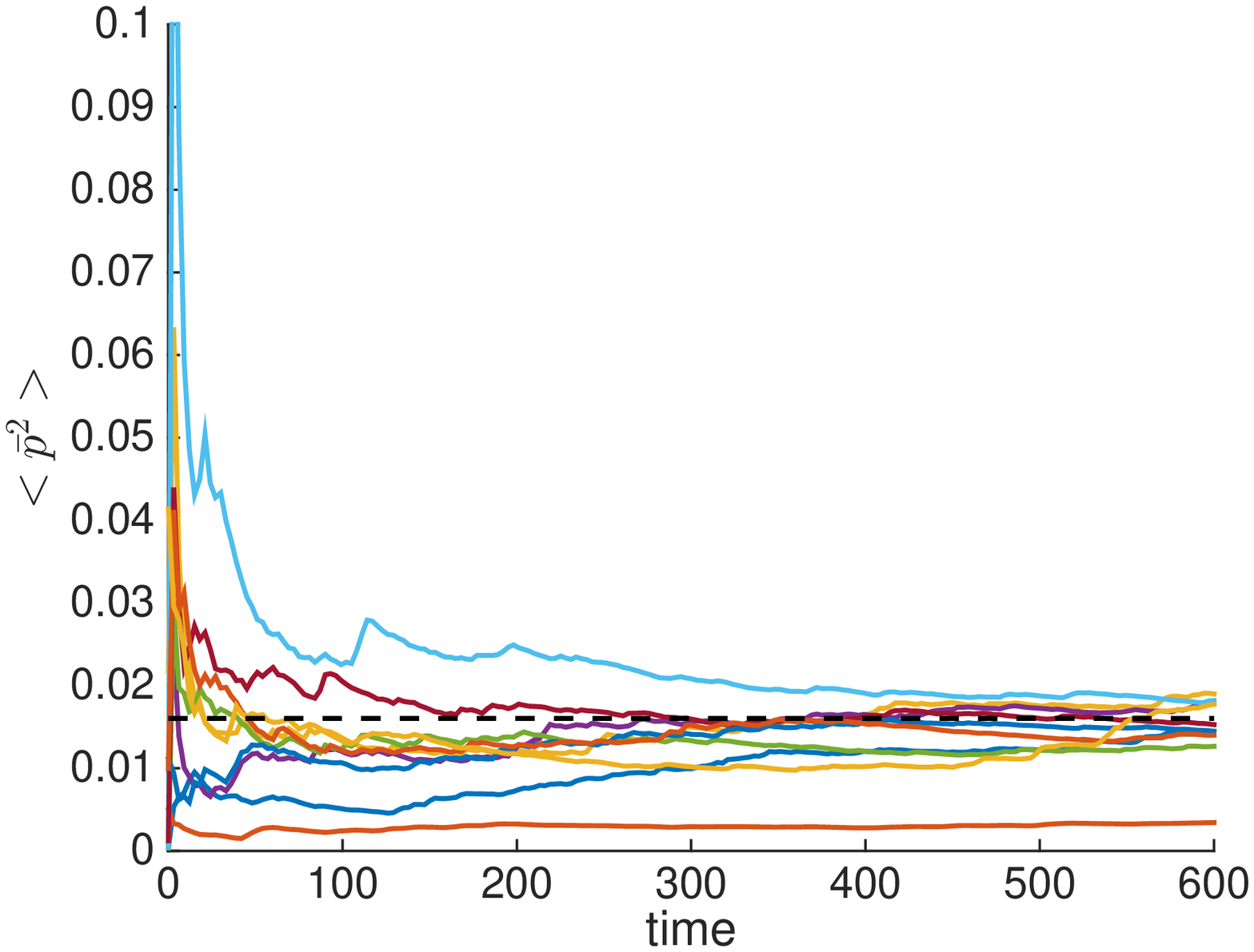} 
   \includegraphics[scale=0.38]{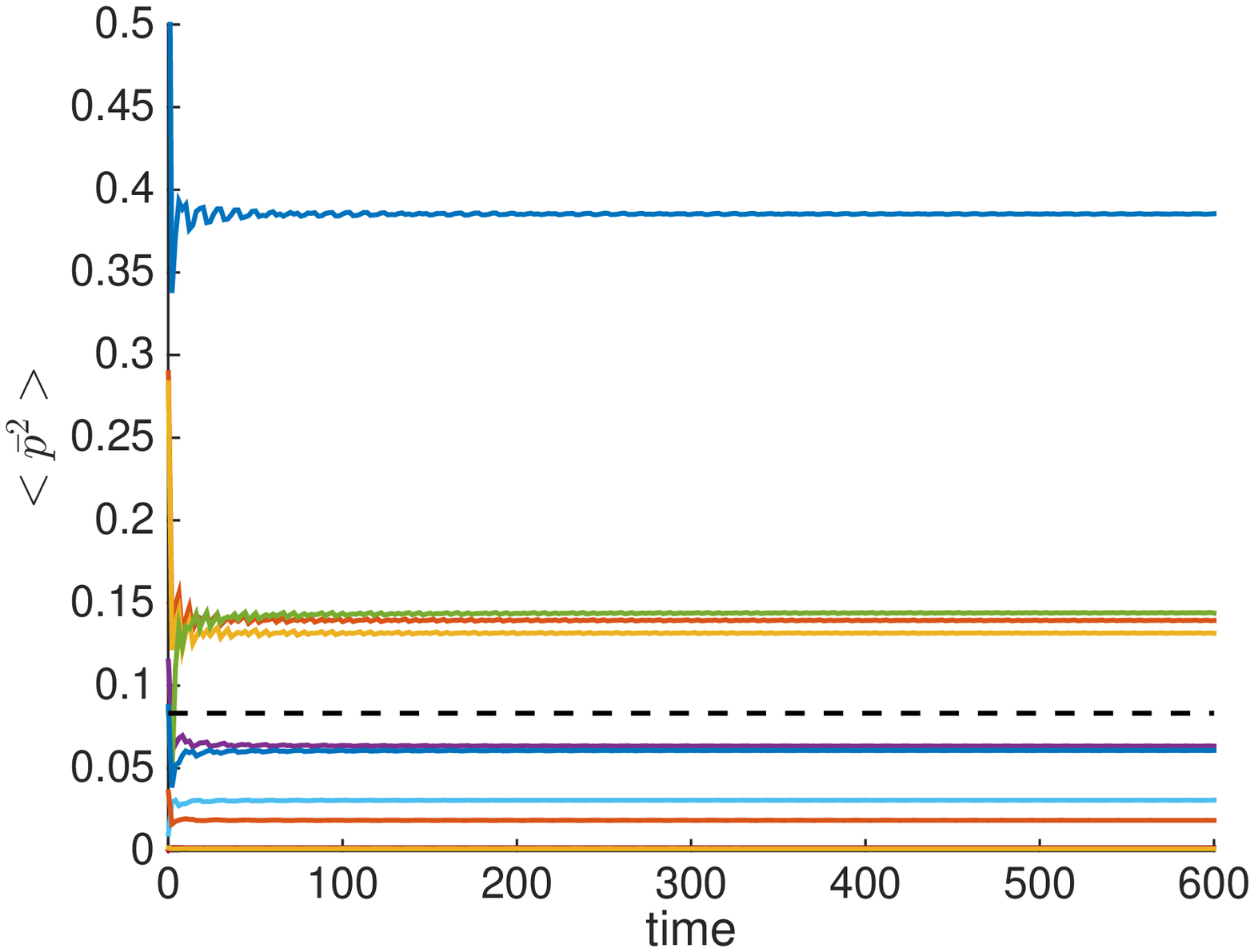}
   \includegraphics[scale=0.38]{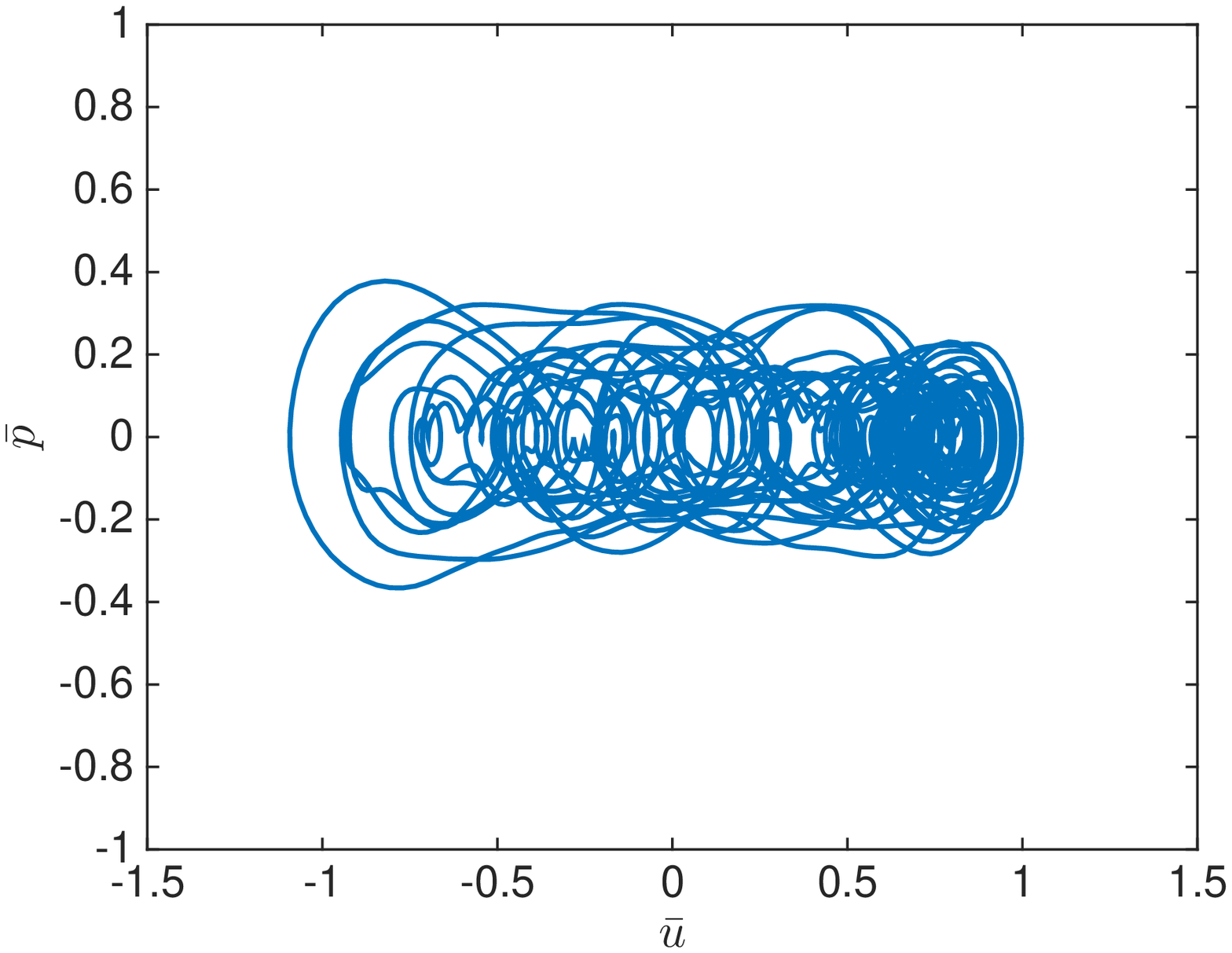} 
   \includegraphics[scale=0.38]{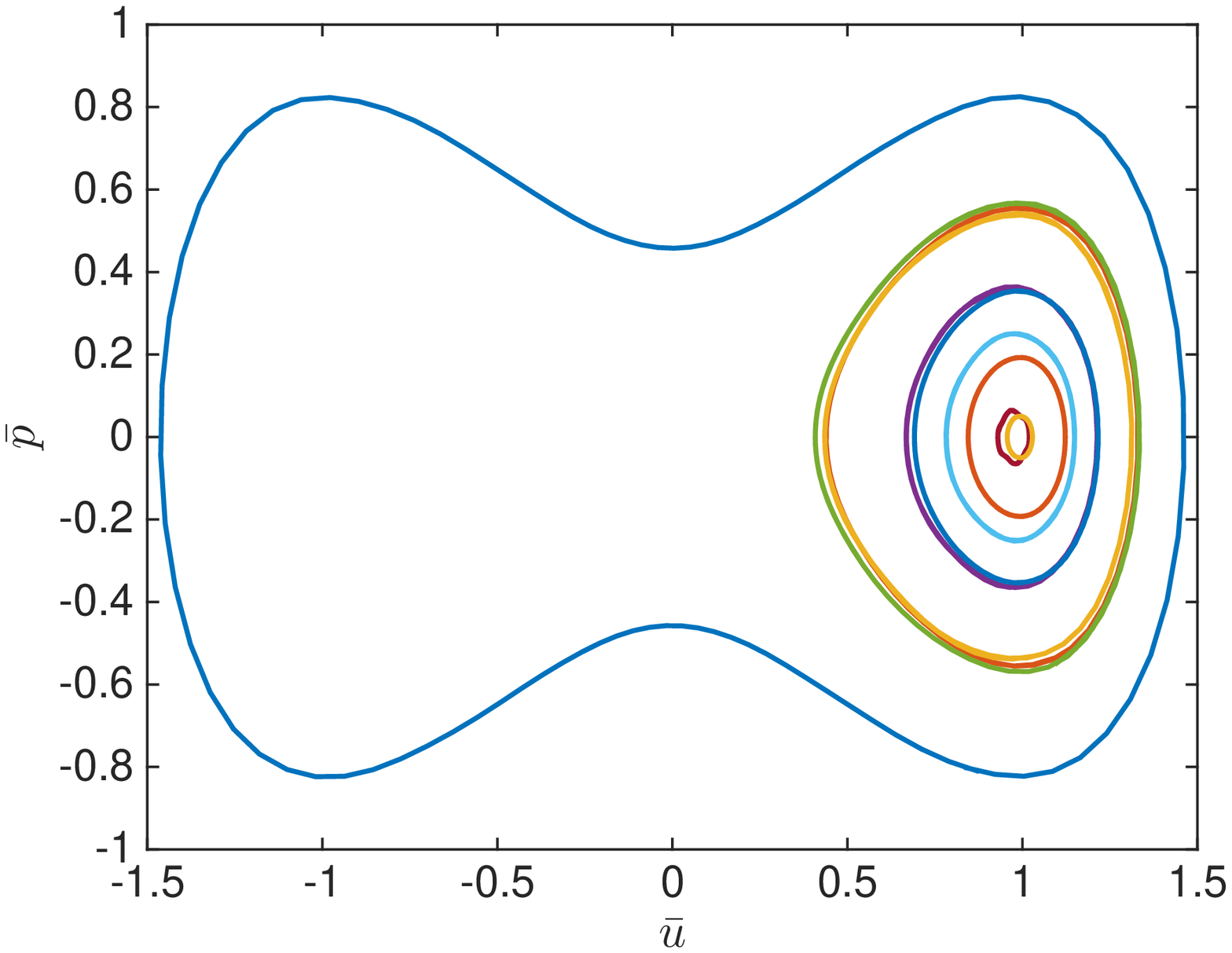}
   \caption{ Results of numerical simulations of the finite dimensional system $\eqref{discrete_sys}$ with $N=1024$.  Simulations starting from 9 different initial conditions, selected according to \eqref{IC}, with $\beta\Delta E = 3$, (left) $\delta = 0.05$ (right) $\delta = 1$, were run until $t=600$ (much larger than the mean transition time).  (top) To emphasize ergodicity in the case of small $\delta$, the time-averaged value of $\bar{p}^2=\sum_{j=1}^N p_j^2/N $ is plotted, and its value of $1/\beta$ obtained from averaging over the invariant measure \eqref{full_IM0} is shown as a dotted line.  (bottom) Dynamics of the same numerical simulations plotted in terms of the reduced variables $\bar{u}=\sum_{j=1}^N u_j^2/N$ and $\bar{p}=\sum_{j=1}^N p_j/N$.  (bottom left) Single trajectory, (bottom right) all 9 trajectories.}
   \label{fig:reduced}
\end{figure}

However, at larger values of $\delta$ we do not see support for ergodicity over the finite simulation time which is orders of magnitude greater than the expected mean transition time.  This point is emphasized in the top row of plots in Fig.~\ref{fig:reduced}, showing that long time averages seem to be independent of initial conditions for small $\delta$ (top left) yet remain strongly dependent on initial conditions for times much longer than an expected mean transition time for large $\delta$ (top right).

This loss of ergodicity can be explained by a fundamental change in the dynamics as $\delta$ increases.  Not only does increasing $\delta$ reduce the size of the fluctuations in the solution $u(x,t)$, but it also makes this solution more uniform in $x$. This is evident from the energy penalty proportional to $\delta^2$ on the $x$-derivatives in the Hamiltonian \eqref{full_ham}, which also lead to the loss of the non-trivial saddle point for $\delta>1/\pi$.  The solutions to \eqref{full_pde0} at each instance in time become approximated by small fluctuations in space about the quantity $\bar{u}(t) = \int_0^1 u(x,t) dx $.


These effects can be quantified by starting from a large $\beta$ expansion for $u(x,t)$, the infinite dimensional version of what was presented in Sec.~\ref{sec:quadratic}.  Taking
\begin{equation}
u(x,t) = u^+(x) + \frac{1}{\sqrt{\beta}} \tilde{u}(x,t)
\end{equation} 
and expanding the Hamiltonian to  quadratic order in $\beta^{-1/2}$, the resulting canonical distribution describes $\tilde{u}$ as a Gaussian process with covariance function 
\begin{equation}\label{Cxy_NBC}
C(x,y) = \left\{ \begin{array}{ll}
\frac{\ds \big( e^{x/\delta} + e^{-x/\delta} \big) \big( e^{(1-y)/\delta} + e^{-(1-y)/\delta} \big)}{\ds 2\delta\left( e^{1/\delta} - e^{-1/\delta}\right)} & x \le y \\
 \frac{ \ds \big( e^{y/\delta} + e^{-y/\delta} \big) \big( e^{(1-x)/\delta} + e^{-(1-x)/\delta} \big)}{ \ds 2\delta \left( e^{1/\delta} - e^{-1/\delta}\right)} & x \ge y
\end{array}\right. ,
\end{equation}
derived in appendix~\ref{sec:Cxy}.  The uniformity of $\tilde{u}(x)$ in $x$ is apparent from the $\delta^{-1}$ exponential dependence indicating the covariance between $\tilde{u}(x)$ at two different points decays slower as $\delta$ increases, resulting in a more uniform description.  More precisely, we can write $u(x)$ selected from the invariant measure at each moment of time approximately as a sum of independent Gaussian random variables, $g_n$,
\begin{equation}\label{ux_series}
u(x) \approx \pm 1 + \frac{1}{\sqrt{\beta}}\left[ g_0 + \sum_{n=1}^\infty \frac{\cos(n\pi x)}{\sqrt{1+\delta^2n^2\pi^2}} g_n \right]
\end{equation}
in which $\cos(n\pi x)$ and $1/(1+\delta^2n^2\pi^2)$ for $n=0,1,2\dots$ are the eigenfunctions and eigenvalues of the covariance $C(x,y)$ in \eqref{Cxy_NBC}.  For small $\delta$, a number of the $n\ge 1$ terms are significant, while for large $\delta$ we see that only the constant zeroth order term is significant.  

Another way to interpret this large $\delta$ behavior is that only the first constant mode of the solution contributes significantly to the dynamics, namely
$$
\bar{u}(t) = \int_0^1 u(x,t) dx
\qquad \text{and} \qquad
\bar{p}(t) = \int_0^1 p(x,t) dx.
$$ 
The dynamics of these reduced variables are 
shown in the bottom row of Fig.~\ref{fig:reduced}.  For small $\delta$ in the bottom left plot, a solution trajectory explores the phase-space, while for large $\delta$ in the bottom right plot,  each trajectory with a different randomly chosen initial condition appears to neraly conserve the reduced Hamiltonian,
$$
\bar{H} = \frac12 \bar{p}^2 + V(\bar{u}),
$$
violating the assumption of ergodicity with respect to the canonical measure with the original Hamiltonian. The conservation of the reduced Hamiltonian in the limit of large $\delta$ can be understood by computing the time evolution of the quantity.  Namely,
$$\begin{aligned}
\frac{d\bar{H}}{dt} &= \bar{p}\int_0^1 \frac{dp}{dt} dx + V'(\bar{u}) \int_0^1 \frac{du}{dt} dx \\
 & = \bar{p}\int_0^1 (\delta^2 u_{xx} - V'(u) )dx + V'(\bar{u})\bar{p} =  O(\delta^{-1})
\end{aligned}$$
in which the last equality holds since $\int_0^1 \delta^2 u_{xx} dx = 0 $ by the Neumann boundary conditions, and $V'(u) = V'(\bar{u}) + O(\delta^{-1})$ when $u(x)$ is nearly uniform by the expansion in \eqref{ux_series}.

\section{Conclusions\label{sec:conclusions}}

We have investigated the dynamics of the nonlinear wave equation~\eqref{full_pde00}, with particular interest in metastability and mean transition times.  We discussed how the finite dimensional microcanonical distribution, with energy scaling extrinsically as $N/\beta$, is equivalent to the canonical distribution as $N\to\infty$.  By revisiting transition state theory (TST) for an ensemble of trajectories chosen with the canonical invariant measure of the finite dimensional Hamiltonian system, we derived the mean frequency at which the solutions cross a dividing surface separating the two minimizers of the Hamiltonian.   Extension to the PDE case required we justify the use of a quadratic approximation of the Hamiltonian in infinite dimensions by showing concentration of the canonical invariant measure for arbitrarily large $N$.   Numerical simulations of the finite dimensional Hamiltonian system supported that the TST time applied for almost all initial conditions only if $\delta$ was small.

This change in behavior from large to small system size does not occur abruptly at a critical value of the parameter $\delta$.  Rather, it appears to be a slow transition as the result of an increase in the mixing time (i.e~the time for the system to thermalize over the measure restricted to a localized region around the closest minimizer) with increasing $\delta$.  At small values of $\delta$, the mixing time appears to be much shorter than the computed mean residency time with TST, suggesting that the metastable dynamics could be approximately coarse-grained onto that of a two-state Markov jump process for the two localized regions around the minimizers.  At intermediate values of $\delta$, the increase in the mixing time causes correlations between successive transitions.  Even though the system could still be ergodic with respect to the canonical invariant measure, the time to observe convergence of long-time averages increases.  In fact, the results in Sec.~\ref{sec:numerical} suggest that this mixing time goes to infinity at least as $O(\delta^{_1})$ as $\delta\to\infty$ and, in particular, eventually becomes much bigger than the residency times in the metastable sets.

The simulated dynamics of the deterministic nonlinear wave equation also provides interesting insight into the low-damping, $\gamma\to 0$, regime of the Langevin equation~\eqref{full_damped}.  The mean transition time computed with large deviation theory (LDT) in \eqref{tauSPDE} is equivalent to the TST time \eqref{tau} in the limit as $\gamma\to 0$.  The implied equivalence of trajectory-wise dynamics is supported only by the numerical simulations of the deterministic system in the small $\delta$ regime, in which the mean transition time computed over a single trajectory is well described by the TST time.  In this regime, both the deterministic equation~\eqref{full_pde0} and its stochastic counterpart~\eqref{full_damped} have equivalent transition times.

This is in contrast to the large $\delta$ regime in which the numerical simulations of the deterministic equation effectively behave like those of a  one-dimensional Hamiltonian system, approximately conserving a reduced energy.  This  illuminates what we might expect from the low-damping regime of the Langevin equation~\eqref{full_damped}.  It should behave  like a one-dimensional noisy-damped system, with values of reduced energy that are often less than the value of the saddle point energy and slowly diffusing relative to the timescale of a typical Hamiltonian orbit.  We would therefore expect the method of averaging~\cite{FWbook}, applied to the reduced variables rather than the full system, to characterize the motion of the reduced energy.  Computing the transition times for the reduced energy as we did in \cite{KanEve2013} would produce a prefactor, inversely proportional to the damping, on the exponential transition time that is much larger than the lower bound computed here with TST.  We hope to pursue this further in future work.

We end with a final remark that the derivation of the TST time cannot be extended to the wave equation in multiple ($D\ge 2$) spatial dimensions.  While the system can still be discretized to form a finite dimensional system, the limit towards the continuous system breaks down.  We cannot justify metastability through the concentration of the invariant measure, and the prefactor containing the eigenvalues in \eqref{tau} fails to converge.  These problems, stemming from the eigenvalues of the Laplacian in multiple spatial dimensions, are reminiscent of known problems with solutions to stochastic partial differential equations like~\eqref{full_damped} when $D\ge2$~\cite{Hairer2014,RNT2012}.

\subsection*{acknowledgements}
We would like to thank Weinan E and Gregor Kova\v ci\v c for useful discussions.  This work was supported in part by the NSF grant DMS-1522767 (E. V.-E.).


\appendix

%
\section{Infinite Dimensional Microcanonical Distribution\label{app:micro_IM}}
%

In this appendix, we prove proposition 1.  We show the limit of the microcanonical distribution  as $N\to\infty$, with energy given by $E_N=N/\beta$, leads to the invariant measure written formally in \eqref{full_IM0} by calculating the characteristic function of both $u(x,t)$ and $p(x,t) = u_t(x,t)$.  We do this exactly for a quadratic Hamiltonian, with $V(u(x)) = \frac{1}{2}f(x) u^2(x)$, and explicitly solve for the covariance function of the continuous stochastic process $u(x,t)$ in three special cases; $u_t(x,t)$ is always white noise. 

We begin with the discrete microcanonical distribution \eqref{micro} with energy $E_N$ for a finite dimensional quadratic Hamiltonian, 
\begin{equation}\label{micro_quad}
d\mu_N = c_N^{-1}\delta\left(E_N - \frac{1}{2N}\vec{p}^T\vec{p} - \frac{1}{2} \vec{u}^TA_N\vec{u}\right) d\vec{p}d\vec{u}
\end{equation}
where the $N\times N$ matrix $A_N$ is such that $\vec{u}^TA_N\vec{u} = \frac{1}{N}\sum_{j=0}^N \delta^2N^2 (u_{j+1}-u_j)^2 + f(x_j)u_j^2$ with $u_0$ and $u_{N+1}$ set by the boundary conditions of \eqref{full_pde00}.  As $N\to\infty$ we have $\vec{u}^T A_N \vec{u}\to \int_0^1 u(x)\big( -\delta^2  \frac{d^2}{dx^2} + f(x) \big) u(x) dx$.  Note the matrix $A_N$ is symmetric in the case of either Dirichlet or Neumann boundary conditions.  

We consider the characteristic function for the $\{u_j\}$ first, finding it converges to
\begin{equation}
\phi_{\vec{u}}(\vec{t} ) \to \exp\left(-\frac{1}{2\beta}\int_0^1 \int_0^1 t(x) C(x,y) t(y) dx dy \right),
\end{equation}
a Gaussian process with covariance $C(x,y)/\beta$.  The covariance function is obtained by solving 
$$
\big( \delta^2  \frac{d^2}{dx^2} + f(x) \big)C(x,y) = \delta(x-y)
$$ 
with appropriate boundary conditions.  A few special cases are given in Sec.~\ref{sec:Cxy}.
Then we follow a similar procedure for the $\{p_j\}$ in Sec.~\ref{sec:p}, finding the characteristic function converges to 
\begin{equation}
\phi_{\vec{p}}(\vec{s}) \to \exp\left( -\frac{1}{2\beta} \int_0^1 s(x)^2 dx\right),
\end{equation}
a Gaussian process with covariance $\delta(x-y)/\beta$ (i.e. white noise).  These are combined to arrive at proposition 1.

\subsection{The $\vec{u}$\label{sec:u}}

The characteristic function for the sets of points $\vec{u}$ selected from the microcanonical distribution \eqref{micro_quad} is defined as
\begin{equation}\label{phi_u}
\phi_{\vec{u}}(\vec{t}) = \frac{\ds \int_{\mathbb{R}^{2N}} e^{\frac{i}{N} \vec{ t}^T \vec{u}} \delta\left(E_N - \frac{1}{2N}\vec{p}^T\vec{p} - \frac{1}{2} \vec{u}^T A_N \vec{u}\right) d\vec{p}d\vec{u}}{
\ds \int_{\mathbb{R}^{2N}} \delta\left(E_N -  \frac{1}{2N}\vec{p}^T\vec{p} - \frac{1}{2} \vec{u}^T A_N \vec{u}\right) d\vec{p}d\vec{u}}.
\end{equation}

The integral in the denominator of \eqref{phi_u} is determined first by switching to the coordinates $\vec{q} = A_N^{1/2}\vec{u}$ and rescaling $\vec{p}$ to be $\vec{p}/\sqrt{N}$.  This results in the integral
$$
\frac{N^{N/2}}{\det A_N^{1/2}}\int_{\mathbb{R}^{2N}} \delta\left(E_N - \frac{1}{2}\vec{p}^T\vec{p} - \frac{1}{2}\vec{q}^T\vec{q} \right) d\vec{p}d\vec{q} 
$$
where $\det A_N^{1/2} = \prod_{j=1}^N \sqrt{\lambda_j}$ and the $\lambda_j$ are the eigenvalues of the matrix $A_N$. Then, switching to spherical  coordinates, $r^2 = \sum_{j=1}^N p_j^2$, we have
$$
\frac{N^{N/2}}{\det A_N^{1/2}} S_{N-1} \int_{\mathbb{R}^N} \int_0^\infty \delta\left(E_N - \frac{1}{2}\vec{q}^T\vec{q} -  \frac{r^2}{2} \right)r^{N-1}drd\vec{q} 
$$
where $S_{N-1} = N\pi^{N/2}/\Gamma(N/2+1)$ is the surface area of the sphere in $N$-dimensional space,
and integrating over $r$ we obtain
$$
\frac{N^{N/2}}{\det A_N^{1/2}} S_{N-1} \int_{\mathbb{R}^N} \left(E_N - \frac{1}{2}\vec{q}^T\vec{q} \right)_+^{
N/2-1}d\vec{q},
$$
where $x_+ = \max(0,x)$. 
Last, we switch to spherical coordinates for the $\vec{q}$,
$$
\frac{N^{N/2}}{\det A_N^{1/2}}  (S_{N-1})^2 \int_0^{\sqrt{2E_N}} \left(E_N - \frac{r^2}{2} \right)^{
N/2-1}r^{N-1}dr 
$$
which is equivalent to
\begin{equation}\label{den}
 \frac{N^{N/2}}{\det A_N^{1/2}} (S_{N-1})^2 E_N^{N-1}  \frac{\sqrt{\pi}\;\Gamma(N/2)}{2^{N/2}\Gamma(N/2+1/2)}.
\end{equation}

For the numerator of \eqref{phi_u} we again switch to the coordinates $\vec{q} = A_N^{1/2}\vec{u}$ and rescale the $\vec{p}$ to $\vec{p}/\sqrt{N}$,  obtaining
$$
\frac{N^{N/2}}{\det A_N^{1/2}} \int_{\mathbb{R}^{2N}} e^{\frac{i}{N} \vec{ t}^T A_N^{-1/2}\vec{q}} \delta\left( E_N - \frac{1}{2}\vec{p}^T\vec{p} - \frac{1}{2}\vec{q}^T\vec{q} \right)d\vec{p}d\vec{q}.
$$
Switching the $\vec{p}$ to spherical coordinates and integrating yields
$$
\frac{N^{N/2}}{\det A_N^{1/2}} S_{N-1} \int_{\mathbb{R}^N} e^{\frac{i}{N} \vec{ t}^T A_N^{-1/2}\vec{q}} \left( E_N - \frac{1}{2}\vec{q}^T\vec{q} \right)_+^{N/2-1} d\vec{q}.
$$
To determine the integral over $\vec{q}$, we write the integral in spherical coordinates over all $\hat{\vec{n}}\in\Omega$, the unit normal to the $N$-dimensional sphere, and all $r$,
\begin{equation}\label{num}
\frac{N^{N/2}}{\det A_N^{1/2}} S_{N-1} \int_0^{\sqrt{2E_N}}  \int_\Omega e^{\frac{i}{N} \vec{ t}^T A_N^{-1/2}\hat{\vec{n}}r} \left( E_N - \frac{r^2}{2}\ \right)^{N/2-1} r^{N-1}d\hat{\vec{n}}dr.
\end{equation}
The appearance of the $\hat{n}$ in the exponent does not allow us to integrate exactly.  Before proceeding, we return to our original expression \eqref{phi_u} by taking the ratio of \eqref{num} to \eqref{den}, resulting in
 \begin{equation}\label{phi5}
\phi_{\vec{u}}(\vec{t} )=  \frac{2^{N/2}\Gamma(N/2+1/2)}{E_N^{N-1} \sqrt{\pi}\;\Gamma(N/2)} \frac{1}{S_{N-1}}  \int_0^{\sqrt{2E_N}}  \int_\Omega e^{\frac{i}{N} \vec{ t}^T A_N^{1/2}\hat{\vec{n}}r} \left( E_N - \frac{r^2}{2}\ \right)^{N/2-1} r^{N-1}d\hat{\vec{n}}dr.
\end{equation}

In \eqref{phi5}, the integral over $\hat{\vec{n}}$ divided by the surface area can be thought of as an expectation with respect to the random variable $\hat{\vec{n}}$, the unit vector chosen uniformly over the surface of the unit sphere in $N$-dimensional space,
$$
 \frac{1}{S_{N-1}} \int_\Omega e^{\frac{i}{N} \vec{ t}^T A_N^{1/2}\hat{\vec{n}}r} d\hat{\vec{n}} = 
\mathbb{E}\left[ e^{\frac{i}{N} \vec{ t}^T A_N^{-1/2}\hat{\vec{n}}r} \right] .
$$
  Alternatively, we can express the random variable $\hat{n}$ as $\vec{g}/|\vec{g}|$ where $\vec{g}$ is an $N$-dimensional vector with components that are independent identically distributed mean zero Gaussian random variables.  Intuitively, for large $N$ this random variable can be approximated by $\vec{g}/\sqrt{N}$ (noticed by Poincar\'e~\cite{Poincare}, see also~\cite{McKean1973}). 
More specifically, in the limit as $N\to\infty$ we have that the expectation over $\hat{n}$ scales as
$$
\mathbb{E}\left[ \exp\left({\frac{i}{N} \vec{ t}^T A_N^{-1/2}\hat{\vec{n}}r}\right) \right] \sim \mathbb{E}\left[ \exp\left(\frac{i}{N^{3/2}} \vec{ t}^T A_N^{-1/2}\vec{g}r\right) \right] 
$$
where we use $\sim$ to indicate that the ratio of the two sides approaches one as $N\to\infty$.  This expectation on the right hand side is the characteristic function of a mean zero multivariate Gaussian; it is known explicitly.  Combining the above results we have \begin{equation}\label{phi5b}
 \frac{1}{S_N} \int_\Omega e^{\frac{i}{N} \vec{ t}^T A_N^{1/2}\hat{\vec{n}}r} d\hat{\vec{n}} \sim 
 \exp\left({-\frac{1}{2}\frac{r^2}{N^3} \vec{t}^TA_N^{-1}\vec{t}}\right)
\end{equation}
in the limit of large $N$.

Then, inserting \eqref{phi5b} into \eqref{phi5} for the integration over $\hat{\vec n}$ we are left with the following integral, 
\begin{equation}\label{int_r}
\int_0^{\sqrt{2E_N}} e^{-\frac{1}{2}\frac{r^2}{N^3} \vec{t}^TA_N^{-1}\vec{t}}  \left( E_N - \frac{r^2}{2}\ \right)^{N/2-1} r^{N-1}dr.
\end{equation}
  We would like to use Laplace's method,
$$
\int_a^b g(x) e^{Nf(x)} dx \sim \sqrt{\frac{2\pi}{N|f''(x_0)|}} g(x_0)e^{Nf(x_0)}
$$
if the integral in \eqref{int_r} has the proper scalings as $N\to\infty$.  We notice the scaling of the exponent in \eqref{int_r} is $r^2/N + O(1)$ as $\vec{t}^TA_N^{-1}\vec{t}/N^2$ converges to a double integral, 
$$
\frac{1}{N^2}   \vec{t}^TA_N^{-1}\vec{t} \to \int_0^1 \int_0^1 t(x) C(x,y) t(y) dx dy
$$
where the function $C(x,y)$ is obtained by inverting the operator $A_\infty$, defined such that $\vec{u}^T A_N \vec{u}\to\int_0^1 u(x) A_\infty(x) u(x) dx $.  This exponential term in \eqref{int_r} remains smaller than $e^N$ as $N\to\infty$ provided the maximum value of $r$ does not grow faster than $O(\sqrt{N})$.

The exponential in $N$ component comes from writing  
$$
\left( E_N - \frac{r^2}{2}\ \right)^{N/2-1} r^{N-1} \sim \exp\left(\frac{N}{2}\log(E_N-\frac{r^2}{2}) + N\log r\right)
$$
to leading order in $N$.  We see that this exponent attains its maximum between $r=0$ and $\sqrt{2E_N}$ when $r=\sqrt{E_N}$.  This corresponds to half of the total energy, $E_N$, contained within the variables $\vec{u}$; the equipartition of energy between the $\vec{p}$ and the $\vec{u}$.   

Now, we define the value of the energy, $E_N = N/\beta$, as this is where the canonical measure energy concentrates, as discussed in the main text.  Even with this energy scaling with $N$, and therefore the maximum value of $r = \sqrt{N/\beta}$, the exponential term in \eqref{int_r} remains smaller than $e^N$.  Applying Laplace's method to the integral in \eqref{int_r}, with energy $E_N = N/\beta$, we obtain
$$\begin{aligned}
\int_0^{\sqrt{2N/\beta}} &e^{-\frac{1}{2}\frac{r^2}{N^3} \vec{t}^TA_N^{-1}\vec{t}}  \left( \frac{N}{\beta} - \frac{r^2}{2}\ \right)^{N/2-1} r^{N-1}dr \\
&\sim \left(\frac{N}{\beta}\right)^{N-1}2^{-N/2+1} \sqrt{\frac{2\pi}{N 4}}  \exp\left({-\frac{1}{2}\frac{1}{\beta N^2} \vec{t}^TA_N^{-1}\vec{t}} \right) 
\end{aligned}$$
as $N\to\infty$.  Combining this result together with \eqref{phi5b} into \eqref{phi5} also with $E_N=N/\beta$, we find after simplifying that
$$
\phi_{\vec{u}}(\vec{t} ) \sim  \frac{\Gamma(N/2+1/2)}{\Gamma(N/2)}  \sqrt{\frac{2}{ N}}  \exp\left({-\frac{1}{2\beta}\frac{1}{N^2} \vec{t}^TA_N^{-1}\vec{t}} \right) .
$$
In the limit as $N\to\infty$ the prefactor in this expression goes to unity (the true value of the integral of $(1-r^2/2)^{N/2-1}r^{N-1}$ and its approximation are asymptotically equivalent), leaving
\begin{equation}
\phi_{\vec{u}}(\vec{t} ) \to \exp\left(-\frac{1}{2\beta}\int_0^1 \int_0^1 t(x) C(x,y) t(y) dx dy \right).
\end{equation}
The stochastic process $u(x)$ has the covariance function $C(x,y)/ \beta$.  The covariance function $C(x,y)$ is found by solving 
\begin{equation}\label{eqCxy}
A_\infty C(x,y) = \delta(x-y)
\end{equation}
with appropriate boundary conditions.  We consider special cases of $A_\infty$ in the next section before proceeding to consider the characteristic function of the $\vec{p}$.

\subsection{Example Covariance Functions\label{sec:Cxy}}

Analytic solutions to \eqref{eqCxy} can be obtained in a few special cases.  Of particular interest is the case of the quadratic function $V(u)=\frac12 u^2$, which can be used to approximate the nonlinear potential function $V(u) = \frac14(1-u^2)^2$ near its minimizers.  

For the case of dynamics governed by \eqref{full_pde0}  with $V(u)=\frac{1}{2}u^2$ and Dirichlet boundary conditions, $u(0,t)=u(1,t)=0$, inverting the operator $A_\infty$ amounts to solving
\begin{equation}
\label{Cxy_DEQ}
\left(-\delta^2  \frac{\partial^2}{\partial x^2} + 1 \right) C(x,y) = \delta(x-y)
\end{equation}
with boundary conditions $C(0,y)=C(1,y)=0$ and continuity conditions at $x=y$.  The solution is 
\begin{equation}\label{Cxy_DBC}
C(x,y) = \left\{ \begin{array}{ll}
\frac{\ds \big( e^{x/\delta} - e^{-x/\delta} \big) \big( e^{(1-y)/\delta} - e^{-(1-y)/\delta} \big)}{\ds 2\delta\left( e^{1/\delta} - e^{-1/\delta}\right)} & x \le y \\
 \frac{ \ds \big( e^{y/\delta} - e^{-y/\delta} \big) \big( e^{(1-x)/\delta} - e^{-(1-x)/\delta} \big)}{ \ds 2\delta \left( e^{1/\delta} - e^{-1/\delta}\right)} & x \ge y
\end{array}\right.
\end{equation}
This covariance is nearly constant along a large portion of $x=y$, indicating the solutions look something like an Ornstein Uhlenbeck process, but conditioned to have $u(0,t)=u(1,t)=0$.

When the boundary conditions are changed to Neumann boundary conditions, the solution to \eqref{Cxy_DEQ} is 
\begin{equation}
\label{Cxy_NBC2}
C(x,y) = \left\{ \begin{array}{ll}
\frac{\ds \big( e^{x/\delta} + e^{-x/\delta} \big) \big( e^{(1-y)/\delta} + e^{-(1-y)/\delta} \big)}{\ds 2\delta\left( e^{1/\delta} - e^{-1/\delta}\right)} & x \le y \\
 \frac{ \ds \big( e^{y/\delta} + e^{-y/\delta} \big) \big( e^{(1-x)/\delta} + e^{-(1-x)/\delta} \big)}{ \ds 2\delta \left( e^{1/\delta} - e^{-1/\delta}\right)} & x \ge y
\end{array}\right.
\end{equation}

To understand what these covariance functions describe, consider no potential function, $V(u)=0$, and Dirichlet boundary conditions.  The covariance is given by the solution to
$$
-\delta^2  \frac{\partial^2}{\partial x^2}   C(x,y) = \delta(x-y)
$$
with boundary conditions $C(0,y)=C(1,y)=0$, which is
\begin{equation}
C(x,y) = \frac{\min(x,y) - xy}{\delta^2},
\end{equation}
the covariance function of a Brownian Bridge on the unit interval.

\subsection{The $\vec{p}$\label{sec:p}}

We can follow the exact same procedure to compute the characteristic function for the $\vec{p}$,
\begin{equation}\label{phi_p}
\phi_{\vec{p}}(\vec{s}) = \frac{\ds \int e^{\frac{i}{N} \vec{ s}^T \vec{p}} \delta\left(E_N - \frac{1}{2N}\vec{p}^T\vec{p} - \frac{1}{2} \vec{u}^T A_N \vec{u}\right) d\vec{p}d\vec{u}}{
\ds \int\delta\left(E_N -  \frac{1}{2N}\vec{p}^T\vec{p} - \frac{1}{2} \vec{u}^T A_N \vec{u}\right) d\vec{p}d\vec{u}}.
\end{equation}
Much of the above derivation follows directly, and we briefly detail the steps.  By rescaling $\vec{p}$ to be $\vec p/\sqrt{N}$ and $\vec{q} = A_N^{1/2}\vec{u}$ we have
$$
\phi_{\vec p} (\vec s)= \frac{\ds \int e^{\frac{i}{\sqrt{N}} \vec{ s}^T \vec{p}} \delta\left(E_N - \frac{1}{2}\vec{p}^T\vec{p} - \frac{1}{2} \vec{q}^T \vec{q}\right) d\vec{p}d\vec{q}}{
\ds \int\delta\left(E_N -  \frac{1}{2}\vec{p}^T\vec{p} - \frac{1}{2} \vec{q}^T \vec{q}\right) d\vec{p}d\vec{q}}.
$$
switching to spherical coordinates and replacing $\hat{\vec{n}}$ with the scaling $\vec{g}/\sqrt{N}$ within the expectation results in 
$$
\phi_{\vec p} (\vec s)= \frac{\ds \int_0^{\sqrt{2E_N}} \mathbb{E}\left[ e^{\frac{i}{\sqrt{N}} \vec{ s}^T \vec{g} r}\right] \left(E_N - \frac{r^2}{2} \right)_+^{
N/2-1} r^{N-1} dr }{
\ds  \int_0^{\sqrt{2E_N}}\left(E_N - \frac{r^2}{2} \right)_+^{
N/2-1} r^{N-1}dr }.
$$
Using that the expectation
$$
\mathbb{E}\left[ \exp(\frac{i}{N} \vec{ s}^T \vec{g}r) \right] =
 \exp(-\frac{r^2}{2N^2} \vec{ s}^T \vec{s})
$$
and using Laplace's method for the integral over $r$ with the scaling of $E_N=N/\beta$ we arrive at
\begin{equation}
\phi_{\vec{p}}(\vec{s}) \to \exp\left( -\frac{1}{2\beta} \int_0^1 s(x)^2 dx\right)
\end{equation}
in the limit as $N\to\infty$. (Note that $\frac{1}{N}\vec{s}^T\vec{s} \to \int_0^1 |s(x)|^2 dx$).  Thus we see that $p(x)$ is white noise with covariance function $\delta(x-y)/ \beta $.

%
\section{Quadratic Approximation\label{app:quad}}
%

In this appendix we prove propositions \ref{prop:concentration} and \ref{prop:asymptotic}.  We show in appendix \ref{app:concentration_c} that an arbitrarily large amount of the invariant measure \eqref{full_IM0} can be contained within an arbitrarily small region around a minimizer of the Hamiltonian for large enough values of $\beta$ by developing a uniform in $N$ concentration of the discrete quadratic canonical measure.  This allows for both the argument that the solutions to \eqref{full_pde0} display metastability, as well as the argument that the potential function $V(u)$ can be replaced with a quadratic approximation in order to complete the integration in determining the time between transitions of the solutions to \eqref{full_pde0}.   Appendix \ref{app:integrals} contains the details of computing the integrals to arrive at Proposition \ref{prop:asymptotic} and Conjecture \ref{prop:asymptotic_micro}, the counterpart to Proposition \ref{prop:asymptotic} for the microcanonical measure found in appendix~\ref{app:micro}.

\subsection{Concentration of the Canonical Measure\label{app:concentration_c}}

In the following, we show the concentration of the marginal canonical measure for $\vec b^\pm= (b^\pm_1, \dots b^\pm_N)$ as stated in Proposition~\ref{prop:concentration} in the main text, and also for the marginal canonical measure restricted to the surface  $S= \{ ( \tilde{\vec p}, \vec a): a_1 = 0\}$ in the eigenvector basis, stated in Proposition~\ref{lemma:c_S} below.

We first point out that the constraint on $\sum_{j=1}^N \delta_j^2$ is to ensure the corner of the box remains a finite distance away from the center of the box regardless of the value of $N$, in particular as $N\to\infty$.  Therefore the measure is being contained within a localized region.
In particular, the sides of the box must decay faster than $1/\sqrt{j}$,
\begin{equation}\label{diag}
\delta_j^2 \sim \frac{C}{j^{1 + \epsilon}},
\end{equation}
for some $\epsilon > 0$ and $C>0$.

We drop the $\pm$ and proceed to calculate
\begin{equation}\label{eq:0}
\frac{\int_{D_\delta^N} e^{-\frac{\beta}{2}\sum_{j=1}^N \lambda_j b_j^2} d \vec b }
{\int_{\mathbb{R}^N} e^{-\frac{\beta}{2}\sum_{j=1}^N \lambda_j b_j^2} d \vec b} .
\end{equation}
by squaring everything so each integral over $b_j$ becomes an integral over $b_j$ and $a_j$:
\begin{equation*}\label{eq:1}
\frac{\ds \int_{-\delta_j}^{\delta_j} \int_{-\delta_j}^{\delta_j} e^{-\frac{\beta}{2}(\lambda_jb_j^2 + \lambda_{j}a_{j}^2)}db_jda_{j}}
{\ds \int_{-\infty}^\infty \int_{-\infty}^\infty e^{-\frac{\beta}{2}(\lambda_jb_j^2 + \lambda_{j}a_{j}^2)}db_jda_{j}}.
\end{equation*}
Then, rescaling by the eigenvalues and noting the  integral over the 2D square with sides $2\delta_j\sqrt{\lambda_j}$ is greater than the integral over the circle with radius $\delta_j\sqrt{\lambda_j}$, gives
\begin{equation*}
\frac{\ds \int_{-\delta_j\sqrt{\lambda_j}}^{\delta_j\sqrt{\lambda_j}} \int_{-\delta_j\sqrt{\lambda_j}}^{\delta_j\sqrt{\lambda_j}} e^{-\frac{\beta}{2}(b_j^2 + a_{j}^2)}db_jda_{j}}
{\ds \int_{-\infty}^\infty \int_{-\infty}^\infty e^{-\frac{\beta}{2}(b_j^2 + a_{j}^2)}db_jda_{j}} 
>
\frac{\ds \int_0^{\delta_j\sqrt{\lambda_j}} e^{-\frac{\beta}{2} r^2}r dr}{\ds \int_0^\infty e^{-\frac{\beta}{2} r^2}r dr} =  \left( 1 - e^{-\frac{\beta}{2}\delta_j^2\lambda_j} \right) 
\end{equation*}
for each $j=1\dots N$.

Combining all the integrals, which are identical, we are left with showing that there exists a $\beta_1$ such that
\begin{equation}\label{eq:3}
 \prod_{j=1}^N \left( 1 - \exp\left(-\frac{\beta}{2}\delta_j^2\lambda_j\right)  \right)  > (1-\hat{\delta})^2
\end{equation}
can be satisfied for all $N>0$ and for all $\beta>\beta_1$.

The product in \eqref{eq:3} converges if the sum
\begin{equation}\label{eq:4}
\sum_{j=1}^N  \log \left( 1 - \exp\left(-\frac{\beta}{2}\delta_j^2\lambda_j\right)  \right) 
\end{equation}
converges as $N\to\infty$.  
We can compare the decay of the terms in this sum to the decay of the terms in a known converging series, $\sum_{j=1}^\infty j^{-p}$ when $p>1$.  We note the terms in \eqref{eq:4} are negative, so we add a negative sign and compare.  We are interested if 
\begin{equation}
\lim_{j\to\infty} \frac{-\log \left( 1 - e^{-\frac{\beta}{2}\delta_j^2 \lambda_j} \right) }{1/j^p} < \infty
\end{equation}
Taking the derivative of top and bottom yields
\begin{equation}\label{eq:5}
\lim_{j\to\infty} \frac{ j^{p+1}\frac{\beta}{2}\delta_j^2\lambda_j e^{-\frac{\beta}{2}\delta_j^2\lambda_j} } {p\left( 1 - e^{-\frac{\beta}{2}\delta_j^2 \lambda_j} \right) } = 0,
\end{equation}
 provided the exponent in the numerator decays to zero, 
\begin{equation}\label{eq:6}
\delta_j^2\lambda_j \to \infty
\end{equation}
as $j\to \infty$.  

We are left to consider if there exits a box with sides $2\delta_j$ that satisfies both \eqref{diag} and \eqref{eq:6}  simultaneously.  This depends on the values of the $\lambda_j$.  For our application, they are the eigenvalues of the operator \eqref{opA}.  They are dominated by the Laplacian term, and therefore for large $j$ scale like 
\begin{equation}\label{eq:7}
\lambda_j \sim j^2.
\end{equation}
From \eqref{diag} we have that $\epsilon > 0$ and from combining \eqref{diag} and \eqref{eq:7} into \eqref{eq:6} we have that $ 1 - \epsilon  > 0 $.  A range of values of $\epsilon$ are possible, for example, $\epsilon = 1/2$.  Thus it is possible to contain the measure within a finite region around a minimizer of the Hamiltonian.  With a choose of $\epsilon$, the value of $\beta_1$ is determined as the one that satisfies \eqref{eq:3} when $N\to\infty$ (the product is decreasing as $N$ increases).

The integration for the TST frequency also requires the following proposition concerning the concentration of the measure
 restricted to the surface $S$, given by $a_1=0$ in the eigenvector coordinates:
%
\begin{proposition}
\label{lemma:c_S}
Consider the marginal canonical measure restricted to the surface $S$ defined in \eqref{eq:20} with the quadratic Hamiltonian approximation \eqref{eq:42}
and an $N-1$-dimensional box with side edge lengths $2\delta_j$, 
$$
\bar{D}_\delta^{N-1} = \{ \vec a : a_1=0,  -\delta_j \le a_j \le \delta_j\;\; \forall \;\;j=2\dots N\},
$$
satisfying $ \sum_{j=2}^N \delta_j^2$ remaining finite as $N\to\infty$. Assume that the eigenvalues $\lambda^s_j$ solving \eqref{eq:22} satisfy
$$
\lambda^s_j \sim Cj^2 \qquad \text{as $j\to\infty$, for some $C>0$}
$$
Then, for every $\hat{\delta}>0$, there exists $\beta_1>0$ and a set of $\{ \delta_j\}_{j=2}^N $ such that for every $\beta>\beta_1$ and every $N>0$
\begin{equation}
\frac{\int_{\bar{D}_\delta^{N-1}} \exp\left(-\frac{\beta}{2}\sum_{j=2}^N \lambda^s_j a_j^2 \right) da_2\dots da_N }
{\int_{\mathbb{R}^{N-1}} \exp\left( -\frac{\beta}{2}\sum_{j=2}^N \lambda^s_j a_j^2\right) da_2\dots da_N} > 1- \hat{\delta}
\end{equation}
\end{proposition}

The proof of Proposition \ref{lemma:c_S} follows exactly, replacing all summations above from $j=2\dots N$.  In particular, the expression in \eqref{eq:4} becomes
\begin{equation}\label{eq:a8}
\sum_{j=2}^N  \log \left( 1 - \exp\left(-\frac{\beta}{2}\delta_j^2\lambda^s_j\right)  \right) 
\end{equation}
and since the eigenvalues 
 $\lambda^s_j$ have the same scaling as \eqref{eq:7}, we arrive at the same conclusion on the convergence of \eqref{eq:a8} and therefore the concentration of the measure.

\subsection{Evaluation of the Integrals\label{app:integrals}}

Having justified the use of a quadratic approximation as $\beta\to\infty$ for the canonical invariant measure, we compute the integration in  \eqref{eq:om2_mc} or \eqref{eq:om2_c} required to compute the TST frequencies.  First, we want to rewrite the integrals in the eigenvector basis.  The normal $\hat{n}(\vec u)\in \mathbb{R}^N$, is given by $(1,0,\dots 0)$  and therefore $\vec p \cdot \hat{n}(\vec u) = \tilde{p}_1$.  
Changing the integration variables to the eigenvector basis is more than just a rotation, as
we defined the eigenvector basis to be normalized to $N$ rather than one.  This results in the integrals in \eqref{eq:om2_mc} and \eqref{eq:om2_c} gaining a factor of $N^{(N-1)/2}$, a factor of $\sqrt{N}$ for each of the $j=2\dots N$ eigenvector coordinates.  Together with the asymptotic in $\beta$ expansion of the measures we have from \eqref{eq:om2_mc}
\begin{equation}\begin{aligned}
\nu_S^m &\sim c_N^{-1} N^{(N-1)/2} \int_{\mathbb{R}^{N-1}} \max(\tilde{p}_1,0)\\
& \times \delta\Big( E_N - E_N^s - \sum_{j=1}^N \frac{\tilde{p}_j^2}{2N} - \sum_{j=2}^N\frac{\lambda_j^sa_j^2}{2} \Big)d\tilde{p}_1\dots d\tilde{p}_N da_2\dots da_N \\
& = c_N^{-1} 2^{N-2} N^{N} \prod_{j=2}^N \frac{1}{\sqrt{\lambda^s_j}} S_{2N-3} \frac{  \left( E_N - E_N^s \right)^{N-1}}{N-1}
\end{aligned}\end{equation}
and from \eqref{eq:om2_c}
\begin{equation}\begin{aligned}
\nu_S^c & \sim C_N^{-1} N^{(N-1)/2}\int_{\mathbb{R}^{N-1}} \max(\tilde{p}_1,0)\\
& \times \exp\Big( - \beta E_N^s - \beta\sum_{j=1}^N \frac{\tilde{p}_j^2}{2N} -\beta \sum_{j=2}^N\frac{\lambda_j^s}{2}a_j^2 \Big)d\tilde{p}_1\dots d\tilde{p}_N da_2\dots da_N \\
& = C_N^{-1}  \frac{N^{N}}{\beta^{N}} \prod_{j=2}^N \frac{1}{\sqrt{\lambda^s_j}} S_{2N-3} 2^{N-2}(N-2)!e^{-\beta E_N^s} 
\end{aligned}\end{equation}
where $S_{n}=2\pi^{(n+1)/2}/\Gamma((n+1)/2)$ is the surface area of the $n$-dimensional sphere (embedded in $n+1$ dimensional space) of unit radius. 

For the normalization constants, the concentration of the measure also allows us to extend the integration to all space with little error, resulting in the asymptotic in $\beta$ expressions
\begin{equation}\begin{aligned}
c_N^\pm &\sim N^{N/2} \int_{\mathbb{R}^{2N}} \delta\Big( E_N - E_N^\pm - \frac{1}{2}\sum_{j=1}^N[ \frac{\tilde{p}_j^2}{N} + \lambda_j^\pm(b_j^\pm)^2] \Big)d\tilde{p}_1\dots d\tilde{p}_N db_1^\pm\dots db_N^\pm 
\end{aligned}\end{equation}
and
\begin{equation}\begin{aligned}
C_N^\pm &\sim N^{N/2}\int_{\mathbb{R}^{2N}} \exp\Big( - \beta E_N^\pm - \frac{\beta}{2} \sum_{j=1}^N[ \frac{\tilde{p}_j^2}{N} + \lambda_j^\pm(b_j^\pm)^2] \Big)d\tilde{p}_1\dots d\tilde{p}_N db_1^\pm\dots db_N^\pm 
\end{aligned}\end{equation}
Notice the factor of $N^{N/2}$ due to the change of the integration to the eigenvector basis.  In our case, due to symmetry, $c_N^+=c_N^-$ and $C_N^+=C_N^-$ ($E^+_N=E^-_N$ and $\lambda_j^+ = \lambda_j^-$ for all $j=1\dots N$), but they could be different if $V(u)$ was not symmetric.  Completing the integration we have
\begin{equation}
c_N^\pm \sim 2^{N-1} N^{N}\prod_{j=1}^N \frac{1}{\sqrt{\lambda^\pm_j}} S_{2N-1} ( E_N - E_N^\pm )^{N-1} 
\end{equation}
and
\begin{equation}
C_N^\pm \sim 2^{N-1} \frac{N^{N}}{\beta^N}S_{2N-1} (N-1)!  \prod_{j=1}^N \frac{1}{ \sqrt{\lambda^\pm_j}} e^{-\beta E_N^\pm} .
\end{equation}
Thus arriving at the expressions in Proposition \ref{prop:asymptotic} and Conjecture \ref{prop:asymptotic_micro}.

%
\section{Expected Residency Time Using the Microcanonical  Measure\label{app:micro}}
%

Here, we compute the average of $\nu_S^T$ over trajectories with initial conditions chosen with respect to the microcanonical \eqref{micro} invariant measure.  As for the case of the canonical invariant measure, the use of a quadratic approximation of the Hamiltonian is justified by the concentration of the measure.
This concentration is not as obvious for the microcanonical measure \eqref{micro}, and will depend upon the dimension being large enough.   Note that the points that minimize the Hamiltonian are in general not on the surface $H=E_N$, however this surface does contain $\vec u^\pm$ with $\vec p$ non-zero satisfying $H(\vec p, \vec u^\pm)=E_N$.  We can therefore think about the concentration of the marginal measure for $\vec u$ near the lowest potential energy point for as was shown for the canonical distribution.  We do not explicitly bound the integration of the microcanonical measure over a localized region, yet we expect the quadratic approximation to be reasonable at large $N$ due to the shown equivalence between the microcanonical and canonical measures as $N\to\infty$.   We present the results assuming this approximation is valid here.  

We proceed with the calculation from Sec.~\ref{sec:TST}, using the microcanonical \eqref{micro} invariant measure.
As a result of \eqref{eq:24}, the ensemble average of the integrand in~\eqref{eq:k} is invariant in time, and after proper interpretation of integration with respect to the Dirac distribution, we arrive at the analogous equation to \eqref{eq:om2_c} for the expected frequency in the case of the microcanonical distribution: 
\begin{equation}
  \label{eq:om2_mc}
  \nu_S^m =   c_N^{-1}\int_{\{H_N(\vec p, \vec u)=E_N\}\cap S}  
  \max\left( \vec p \cdot \hat n(\vec u) ,0\right) \frac{d\hat{\sigma}(\vec p, \vec u)}{|\nabla H_N(\vec p, \vec u)|} 
\end{equation}
where $d\hat{\sigma}(\vec p, \vec u)$ is the surface element on the intersection of $\{H_N(\vec p, \vec u)=E_N\}$ and $S$, and $\hat n(\vec u)\in \RR^N$ is the unit normal to $\{\vec u : q(\vec u)=0\}$. 

For the normalization constant, we first separate the integration into the two regions, $B_+=\{ (\vec p,\vec u): (\vec u-\vec u^s)\cdot \vec\phi^{(1)} >0\}$ and $B_-=\{ (\vec p,\vec u): (\vec u-\vec u^s)\cdot \vec\phi^{(1)} <0\}$ that partition space, separated by $S$, 
\begin{equation}
\label{eq:43b}
c_N = c_N^+ + c_N^-,
\end{equation}
where
\begin{equation}
\label{eq:44a}
c_N^\pm = \int_{H_N(\vec p, \vec u)=E_N\cap B_\pm} \frac{d\sigma(\vec{p},\vec{u})}{|\nabla H_N(\vec{p},\vec{u})|}  .
\end{equation}

The details of obtaining the following expressions appear in appendix \ref{app:integrals}:
\begin{conjecture}
  \label{prop:asymptotic_micro}  
For the microcanonical ensemble with energy $E_N = N/\beta$,  the asymptotic in $\beta\to\infty$, valid for large $N$, expansion of the integral in \eqref{eq:om2_mc}  is
\begin{equation}
\int_{\{H_N(\vec p, \vec u)=E_N\}\cap S}  
  \max\left( \vec p \cdot \hat n(\vec u) ,0\right) \frac{d\hat{\sigma}(\vec p, \vec u)}{|\nabla H_N(\vec p, \vec u)|}  \sim \frac{2^{N-2} N^{N} S_{2N-3}}{(N-1) \prod_{j=2}^N \sqrt{\lambda^s_j}}    \left( \frac{N}{\beta} - E_N^s \right)^{N-1}
\end{equation}
and the expansion of the integrals in \eqref{eq:44a} are
\begin{equation}
c_N^\pm \sim \frac{2^{N-1} N^{N}S_{2N-1}}{\prod_{j=1}^N \sqrt{\lambda^\pm_j}} \left( \frac{N}{\beta} - E_N^\pm \right)^{N-1} .
\end{equation}
In the above, the eigenvalues $\lambda_j^s$, $j=2\dots N$ solve \eqref{eq:22}, the eigenvalues $\lambda_j^\pm$, solve \eqref{eq:46}
and $S_{n}=2\pi^{(n+1)/2}/\Gamma((n+1)/2)$ is the surface area of the $n$-dimensional sphere (embedded in $n+1$ dimensional space) of unit radius.
\end{conjecture}

Using the above asymptotic in $\beta$ expansions, the TST frequency from \eqref{eq:om2_mc} for the microcanonical ensemble with energy $N/\beta$ is 
\begin{equation}
\nu_S^m \sim \frac{1}{2\pi } \frac{\prod_{j=2}^N \frac{1}{\sqrt{\lambda^s_j}}( 1 - \frac{\beta E_N^s}{N} )^{N-1}}{\prod_{j=1}^N\frac{1}{ \sqrt{\lambda^+_j}}( 1 - \frac{\beta E_N^+}{N} )^{N-1} + \prod_{j=1}^N \frac{1}{ \sqrt{\lambda^-_j}}( 1 - \frac{\beta E_N^-}{N} )^{N-1}} 
\end{equation}
and the mean residency time given by
\begin{equation}
\tau_\pm^m = \frac{c_N^\pm}{c_N \nu_S^m}
\end{equation}
 is
 \begin{equation}\label{tau_finite2}
\tau_\pm^m  \sim 2 \pi \frac{\prod_{j=2}^N \sqrt{\lambda_j^s}}{\prod_{j=1}^N \sqrt{\lambda_j^\pm}}\frac{(1-\frac{\beta E_N^\pm}{N})^{N-1}}{(1-\frac{\beta E_N^s}{N})^{N-1}}.
\end{equation} 
Notice that taking the $N\to\infty$ limit of \eqref{tau_finite2} also converges to \eqref{tau}.
This relies on the fact that the fraction $ (1-\beta E_N^\pm/N)^{N-1}/(1-\beta E_N^s/N )^{N-1}$ converges to $e^{\beta\Delta E^\pm}$, defining the energy barrier height $\mathit{\Delta} E_\pm = U(u^s)-U(u^\pm)$.


\end{document}